\documentclass[leqno,11pt]{amsart}
\usepackage{amsmath,amscd}
\usepackage{epsfig,graphics}
\usepackage{amssymb,latexsym}
\usepackage{mathrsfs}
\usepackage{eucal}

\newtheorem{thm}{\bf Theorem}[section]

\newtheorem{prop}[thm]{\bf Proposition}
\newtheorem{cor}[thm]{\bf Corollary}
\newtheorem{lem}[thm]{\bf Lemma}
\newtheorem{rem}[thm]{\bf Remark}
\newtheorem{ex}[thm]{\bf Example}
\newtheorem{conj}[thm]{\bf Conjecture}

\newcommand{\A}{\mathcal{A}}

\newcommand{\B}{\mathbf{B}}
\newcommand{\cB}{\mathcal{B}}

\newcommand{\cP}{\mathscr{P}}

\newcommand{\pf}{\noindent{\bfseries Proof. }}

\newcommand{\M}{{\mathbf{M}}}
\newcommand{\E}{{\mathcal{E}}}

\newcommand{\cM}{\mathcal{M}}
\newcommand{\F}{\mathcal{F}}

\newcommand{\gl}{\mathfrak{gl}}
\newcommand{\Z}{\mathbb{Z}}
\newcommand{\C}{\mathbb{C}}
\newcommand{\h}{\mathfrak{h}}
\newcommand{\te}{\widetilde{e}}
\newcommand{\tf}{\widetilde{f}}

\newcommand{\tE}{\widetilde{E}}
\newcommand{\tF}{\widetilde{F}}

\numberwithin{equation}{section}

\begin{document}
\title[ ]
{Crystal bases of modified quantized enveloping algebras and a
double RSK correspondence}
\author{JAE-HOON KWON}
\address{Department of Mathematics \\ University of Seoul   \\  Seoul 130-743, Korea }
\email{jhkwon@uos.ac.kr }


\begin{abstract}
The crystal base of the modified quantized enveloping algebras of
type $A_{+\infty}$ or $A_\infty$ is realized as a set of integral
bimatrices. It is obtained by describing the decomposition of the
tensor product of a highest weight crystal and a lowest weight
crystal into extremal weight crystals, and taking its limit using a
tableaux model of extremal weight crystals. This realization induces
in a purely combinatorial way a bicrystal structure of the crystal
base of the modified quantized enveloping algebras and hence its
Peter-Weyl type decomposition generalizing the classical RSK
correspondence.
\end{abstract}

\maketitle

\section{Introduction}
Let $U_q(\frak{g})$ be the quantized enveloping algebra associated
with a symmetrizable Kac-Moody algebra $\frak{g}$.  In \cite{Lu92},
Lusztig introduced the modified quantized enveloping algebra
$\widetilde{U}_q(\frak{g})=\bigoplus_{\Lambda}U_q(\frak{g})a_\Lambda$,
where $\Lambda$ runs over all integral weight for $\mathfrak{g}$,
and proved the existence of its global crystal base or canonical
basis. In \cite{Kas94'}, Kashiwara studied the crystal structure of
$\widetilde{U}_q(\frak{g})$ in detail, and showed that
\begin{equation*}\label{Balambda}
\B(U_q(\frak{g})a_\Lambda)\simeq \B(\infty)\otimes T_\Lambda\otimes \B(-\infty),
\end{equation*}
where $\B(U_q(\frak{g})a_\Lambda)$ denotes the crystal base of
$U_q(\frak{g})a_\Lambda$, $\B(\pm\infty)$  is the crystal base of
the negative (resp. positive) part of $U_q(\frak{g})$ and
$T_\Lambda$ is a crystal with a single element $t_\Lambda$ with
$\varepsilon_i(t_\Lambda)=\varphi_i(t_\Lambda)=-\infty$. It is also
shown that the anti-involution $\ast$ on $\widetilde{U}_q(\frak{g})$
provides the crystal
$\B(\widetilde{U}_q(\frak{g}))=\bigsqcup_\Lambda\B(\infty)\otimes
T_\Lambda\otimes \B(-\infty)$ with another crystal structure called
$\ast$-crystal structure and therefore a regular
$(\frak{g},\frak{g})$-bicrystal structure \cite{Kas94'}.  With
respect to this bicrystal structure, a Peter-Weyl type decomposition
for $\B(\widetilde{U}_q(\frak{g}))$ was obtained when it is of
finite type or affine type at non-zero levels by Kashiwara
\cite{Kas94'} and of affine type at level zero by Beck and Nakajima
\cite{BN} (see also \cite{Na97,Na} for partial results). It is also
shown in \cite{Kas94'} that the crystal graph of the quantized
coordinate ring for $\mathfrak{g}$ \cite{Kas93} is a subcrystal of
$\B(\widetilde{U}_q(\frak{g}))$, and equal to
$\B(\widetilde{U}_q(\frak{g}))$ if and only if $\mathfrak{g}$ is of
finite type.

The main purpose of this work is to study the crystal structure of
the modified quantized enveloping algebras in type $A$ of infinite
rank in terms of Young tableaux and understand its connection with
the classical RSK correspondence. Note that the essential ingredient
for understanding the structure of $\widetilde{U}_q(\frak{g})$ is
the notion of extremal weight $U_q(\mathfrak{g})$-module introduced
by Kashiwara \cite{Kas94'}. An extremal weight module associated
with an integral weight $\lambda$ is an integrable
$U_q(\mathfrak{g})$-module, which is a generalization of a highest
weight and a lowest weight module, and it also has a (global)
crystal base. Our approach is based on the combinatorial models for
extremal weight crystals of type $A_{+\infty}$ and $A_\infty$
developed in \cite{K09,K09-2}.

From now on, we denote $\mathfrak{g}$ by $\gl_{>0}$ and
$\gl_{\infty}$ when it is a general linear Lie algebra of type
$A_{+\infty}$ and $A_\infty$, respectively. The main result in this
paper gives an explicit combinatorial realization  of
$\B(\infty)\otimes T_\Lambda\otimes \B(-\infty)$ for all integral
$\gl_{>0}$-weights and all level zero integral
$\gl_{\infty}$-weights $\Lambda$ as a set of certain bimatrices,
which implies directly Peter-Weyl type decompositions of
$\B(\widetilde{U}_q(\frak{gl}_{>0}))$ and the level zero part of
$\B(\widetilde{U}_q(\frak{gl}_{\infty}))$ without using the
$\ast$-crystal structure.

Let us state our results more precisely. Let $\cM$ be the set of
$\mathbb{N}\times \mathbb{N}$ matrices with finitely many
non-negative integral entries. Recall that $\cM$ has a
$\gl_{>0}$-crystal structure where each row of a matrix in $\cM$ is
identified with a single row Young tableau or a crystal element
associated  with the symmetric power of the natural representation
(cf.\cite{DK}). Let $\cM^\vee=\{\,M^\vee\,|\,M\in\cM\,\}$ be the
dual crystal of $\cM$. For each integral weight $\Lambda$, let
\begin{equation*}
\widetilde{\cM}_\Lambda=\{\,M^\vee\otimes N
\,|\,{\rm wt}(N^t) - {\rm wt}(M^t) =\Lambda \,\}\subset \cM^\vee\otimes\cM.
\end{equation*}
Here ${\rm wt}$ denotes the weight with respect to $\gl_{>0}$-crystal structure
and $A^t$ denotes the transpose of $A\in\cM$. Then our main result (Theorem \ref{Momega}) is
\begin{equation*}\label{main-1}
\widetilde{\cM}_\Lambda \simeq \B(\infty)\otimes T_{\Lambda}\otimes
\B(-\infty).
\end{equation*}
Note that a connected component of $\B(\infty)\otimes
T_\Lambda\otimes \B(-\infty)$  is in general an extremal weight
$\gl_{>0}$-crystal, and  an extremal weight $\gl_{>0}$-crystal is
isomorphic to the tensor product of a lowest weight crystal and a
highest weight crystal \cite{K09}. The crucial step in the proof is
the description of the tensor product
$\B(\Lambda')\otimes\B(-\Lambda'')$ for dominant integral weights
$\Lambda'$, $\Lambda''$ with $\Lambda=\Lambda'-\Lambda''$ in terms
of skew Young bitableaux (Proposition \ref{isomorphism-1+infty}),
and its embedding into $\B(\Lambda'+\xi)\otimes\B(-\xi-\Lambda'')$
for arbitrary integral dominant weight $\xi$  (Proposition
\ref{isomorphism-2+infty}). In fact,
$\B(\Lambda'+\xi)\otimes\B(-\xi-\Lambda'')$ is realized as a set of
skew Young bitableaux whose shapes are almost horizontal strips as
$\xi$ goes to infinity. This establishes the above isomorphism and
as a consequence
\begin{equation*}\label{main-2}
\B(\widetilde{U}_q(\gl_{>0})) \simeq \cM^\vee\otimes \cM,
\end{equation*}
since $\bigsqcup_{\Lambda}\widetilde{\cM}_\Lambda=\cM^\vee\otimes \cM$.

Now, for partitions $\mu,\nu$, let $\cB_{\mu,\nu}$ be the extremal
weight crystal with the Weyl group orbit of its extremal weight
corresponding to the pair $(\mu,\nu)$. Note that
$\cB_{\mu,\emptyset}$ (resp. $\cB_{\emptyset,\nu}$) is a highest
(resp. lowest) weight crystal and it is shown in \cite{K09} that
$\cB_{\mu,\nu}\simeq \cB_{\emptyset,\nu}\otimes
\cB_{\mu,\emptyset}$. Then a $(\gl_{>0},\gl_{>0})$-bicrystal
structure of $\cM$ and $\cM^\vee$ arising from the RSK
correspondence naturally induces a $(\gl_{>0},\gl_{>0})$-bicrystal
structure of $\B(\widetilde{U}_q(\gl_{>0}))$ and  the following
Peter-Weyl type decomposition (Corollary \ref{PeterWeyl+infty})
\begin{equation*}\label{main-3}
\B(\widetilde{U}_q(\gl_{>0})) \simeq
\bigsqcup_{\mu,\nu}\cB_{\mu,\nu}\times\cB_{\mu,\nu}.
\end{equation*}
Hence the decomposition of $\B(\widetilde{U}_q(\gl_{>0}))$ into
extremal weight crystals can be understood as the tensor product of
two RSK correspondences, which are dual to each other as a
$(\gl_{>0},\gl_{>0})$-bicrystal.

Next, we prove analogues for the level zero part of
$\B(\widetilde{U}_q(\mathfrak{gl}_\infty))$. This is done by taking
the limit of the results in $\gl_{>0}$. In this case, $\cM$ is
replaced by $\Z\times \Z$-matrices and $\cB_{\mu,\nu}$ is replaced
by the level zero extremal weight $\gl_\infty$-crystal with the same
parameter $(\mu,\nu)$. Finally, we conjecture that the second
crystal structures arising from the RSK correspondence is compatible
with the dual of $\ast$-crystal structure.

There are several nice combinatorial descriptions of $\B(\infty)$
for $\gl_{>0}$ and $\gl_{\infty}$,  by which we can understand the
structure of the modified quantized enveloping algebras (see e.g.
\cite{HN,LTV,NZ,R}). But our combinatorial description of
$\B(U_q(\frak{gl}_{>0})a_\Lambda)$ and
$\B(U_q(\frak{gl}_{\infty})a_\Lambda)$ explains more directly its
connected components, projections onto tensor products of a highest
weight crystal and a lowest weight crystal, and a bicrystal
structure on $\B(U_q(\frak{gl}_{>0}))$ and
$\B(U_q(\frak{gl}_{\infty}))$.

The paper is organized as follows. In Section 2, we give necessary
background on crystals.  In Section 3, we recall some combinatorics
of Littlewood-Richardson tableaux from a view point of crystals,
which is necessary for our later arguments. In Section 4, we review
a combinatorial model of extremal weight $\gl_{>0}$-crystals and
their non-commutative Littlewood-Richardson rules, and then in
Section 5  we prove the main theorem. In Section 6, we recall the
combinatorial model for extremal weight $\gl_\infty$-crystals and
describe the Littlewood-Richardson rule for the tensor product of a
highest weight crystal and a lowest weight crystal into extremal
ones. In Section 7, we prove analogues for the level zero part of
$\B(\widetilde{U}_q(\mathfrak{gl}_\infty))$. We remark that the
Littlewood-Richardson rule in Section 6 is not necessary for Section
7, but is of independent interest, which completes the discussion on
tensor product of extremal weight $\gl_\infty$-crystals in
\cite{K09-2}.

\section{Crystals}
\subsection{}

Let $\gl_\infty$ be the Lie algebra of complex matrices
$(a_{ij})_{i,j\in \Z}$ with finitely many non-zero entries, which is
spanned by $E_{ij}$ ($i,j\in\Z$), the elementary matrix with $1$ at
the $i$-th row and the $j$-th column and zero elsewhere.

Let $\h=\bigoplus_{i\in \Z}\C E_{ii}$ be the Cartan subalgebra of
$\gl_\infty$ and $\langle\cdot,\cdot\rangle$ denote the natural
pairing on $\frak{h}^*\times\frak{h}$. We denote by $\{\,
h_i=E_{ii}-E_{i+1,i+1}\,|\, i\in\Z\, \}$ the set of simple coroots,
and denote by $\{\, \alpha_i=\epsilon_i-\epsilon_{i+1} \,|\,
i\in\Z\, \}$ the set of simple roots, where $\epsilon_i\in\h^*$ is
given by $\langle \epsilon_i,E_{jj}\rangle=\delta_{ij}$.

Let $P=\Z\Lambda_0\oplus
\bigoplus_{i\in\Z}\Z\epsilon_i=\bigoplus_{i\in\Z}\Z\Lambda_i$ be the
weight lattice of $\gl_{\infty}$, where $\Lambda_0$ is given by
$\langle\Lambda_0,E_{-j+1,-j+1}\rangle=-\langle\Lambda_0,E_{jj}\rangle=\frac{1}{2}$
($j\geq 1$), and $\Lambda_i=\Lambda_0+\sum_{k=1}^i\epsilon_k$,
$\Lambda_{-i}=\Lambda_0-\sum_{k=i}^{0}\epsilon_k$ for $i\geq 1$. We
call $\Lambda_i$ the $i$-th fundamental weight.

For $k\in\Z$, let $P_k=k\Lambda_0+\bigoplus_{i\in\Z}\Z\epsilon_i$ be
the set of integral weights of level $k$. Let $P^+=\{\,\Lambda\in
P\,|\,\langle\Lambda, h_i\rangle\geq 0,\
i\in\Z\,\}=\sum_{i\in\Z}\Z_{\geq 0}\Lambda_i$ be the set of dominant
integral weights. We put $P^+_k=P^+\cap P_k$ for $k\in \Z$. For
$\Lambda=\sum_{i\in\Z}c_i\Lambda_i\in P$, the level of $\Lambda$ is
$\sum_{i\in\Z}c_i$.  If we put $\Lambda_\pm=\sum_{i; c_i\gtrless
0}c_i\Lambda_i$, then $\Lambda=\Lambda_+-\Lambda_-$ with
$\Lambda_\pm\in P^+$.

For $i\in\Z$, let $r_i$ be the simple reflection given by
$r_i(\lambda)=\lambda-\langle \lambda,h_i \rangle\alpha_i$ for
$\lambda\in \frak{h}^*$. Let $W$ be the Weyl group of
$\gl_{\infty}$, that is, the subgroup of $GL(\frak{h}^*)$ generated
by $r_i$ for $i\in\Z$.

For $p,q\in\Z$, let $[p,q]=\{\,p,p+1,\ldots,q\,\}$ ($p<q$), and
$[p,\infty)=\{\,p,p+1,\ldots\}$,
$(-\infty,q\,]=\{\,\ldots,q-1,q\,\}$. For simplicity, we denote
$[1,n]$ by $[n]$ ($n>1$). For an interval $S$ in $\Z$,  let
$\gl_{S}$ be the subalgebra of $\gl_{\infty}$ spanned by $E_{ij}$
for $i,j\in S$.  We denote by ${S}^\circ$ the index set of simple
roots for $\gl_{S}$. For example, $[p,q]^{\circ}=\{\,p,\ldots,
q-1\,\}$. We also put $\gl_{>r}=\gl_{[r+1,\infty)}$ and
$\gl_{<r}=\gl_{(-\infty,r-1]}$ for $r \in \Z$.

\subsection{}\label{crystal}
Let us briefly recall the notion of crystals (see \cite{Kas94} for a
general review and references therein).

Let $S$ be an interval in $\Z$. A {\it $\gl_{S}$-crystal} is a set $B$ together with the maps
${\rm wt}  : B \rightarrow P$, $\varepsilon_i, \varphi_i: B
\rightarrow \mathbb{Z}\cup \{-\infty\}$ and $\te_i, \tf_i: B
\rightarrow B\cup\{{\bf 0}\}$ ($i\in {S}^\circ$) such that for $b\in B$
\begin{itemize}
\item[(1)] $\varphi_i(b) =\langle {\rm wt}(b),h_i \rangle +
\varepsilon_i(b),$

\item[(2)]  $\varepsilon_i(\te_i b) = \varepsilon_i(b) - 1$, $\varphi_i(\te_i b) =
\varphi_i(b) + 1$, ${\rm wt}(\te_ib)={\rm wt}(b)+\alpha_i$ if $\te_i b \neq {\bf 0}$,

\item[(3)] $\varepsilon_i(\tf_i b) = \varepsilon_i(b) + 1$, $\varphi_i(\tf_i b) =
\varphi_i(b) - 1$, ${\rm wt}({\tf_i}b)={\rm wt}(b)-\alpha_i$ if $\tf_i b \neq {\bf 0}$,

\item[(4)] $\tf_i b = b'$ if and only if $b = \te_i
b'$ for $b, b' \in B$,

\item[(5)] $\te_ib=\tf_ib={\bf 0}$ if $\varphi_i(b)=-\infty$,
\end{itemize}
where ${\bf 0}$ is a formal symbol and $-\infty$ is the smallest
element in $\mathbb{Z}\cup\{-\infty\}$ such that $-\infty+n=-\infty$
for all $n\in\mathbb{Z}$. Note that $B$ is equipped with a colored
oriented graph structure, where $b\stackrel{i}{\rightarrow}b'$ if
and only if $b'=\tf_{i}b$ for $b,b'\in B$ and $i\in S^\circ$. We
call $B$  {\it connected} if it is connected as a graph. We call $B$
{\it regular} if $\varepsilon_i(b)={\rm max}\{\,k\,|\,\te_i^kb\neq
{\bf 0}\,\}$ and $\varphi_i(b)={\rm max}\{\,k\,|\,\tf_i^kb\neq {\bf
0}\,\}$ for $b\in B$ and $i\in {S}^\circ$. The {\it dual crystal
$B^\vee$ of $B$} is defined to be the set $\{\,b^\vee\,|\,b\in
B\,\}$ with ${\rm wt}(b^\vee)=-{\rm wt}(b)$,
$\varepsilon_i(b^\vee)=\varphi_i(b)$,
$\varphi_i(b^\vee)=\varepsilon_i(b)$,
$\widetilde{e}_i(b^\vee)=\left(\widetilde{f}_i b \right)^\vee$ and
$\widetilde{f}_i(b^\vee)=\left(\widetilde{e}_i b \right)^\vee$ for
$b\in B$ and $i\in {S}^\circ$. We assume that ${\bf 0}^\vee={\bf
0}$.

Let $B_1$ and $B_2$ be crystals. A {\it morphism}
$\psi : B_1 \rightarrow B_2$ is a map from $B_1\cup\{{\bf 0}\}$ to
$B_2\cup\{{\bf 0}\}$ such that for $b\in B_1$ and $i\in S^\circ$
\begin{itemize}
\item[(1)] $\psi(\bf{0})=\bf{0}$,

\item[(2)] ${\rm wt}(\psi(b))={\rm wt}(b)$,
$\varepsilon_i(\psi(b))=\varepsilon_i(b)$, and
$\varphi_i(\psi(b))=\varphi_i(b)$ if $\psi(b)\neq \bf{0}$,

\item[(3)] $\psi(\te_i b)=\te_i\psi(b)$ if $\psi(b)\neq \bf{0}$ and
$\psi(\te_i b)\neq \bf{0}$,

\item[(4)] $\psi(\tf_i
b)=\tf_i\psi(b)$ if $\psi(b)\neq \bf{0}$ and
$\psi(\tf_i b)\neq \bf{0}$.
\end{itemize}
We call $\psi$ an {\it embedding} and $B_1$ a {\it subcrystal of}
$B_2$ when $\psi$ is injective, and call $\psi$ {\it strict} if
$\psi : B_1\cup\{{\bf 0}\} \rightarrow B_2\cup\{{\bf 0}\}$ commutes
with $\te_i$ and $\tf_i$ for $i\in S^\circ$, where we assume that
$\te_i{\bf 0}=\tf_i{\bf 0}={\bf 0}$. If $\psi$ is a strict
embedding, then $B_2$ is isomorphic to $B_1 \sqcup (B_2\setminus
B_1)$. Note that an embedding between regular crystals is always
strict. For $b_i\in B_i$ ($i=1,2$), we say that {\it $b_1$ is {\rm
($\gl_{{S}}$-)}equivalent to $b_2$}, and write $b_1 \equiv b_2$ if
there exists an isomorphism of crystals $C(b_1)\rightarrow C(b_2)$
sending $b_1$ to $b_2$, where $C(b_i)$ denote the connected
component of $B_i$ including $b_i$.

A {\it tensor product} of $B_1$ and $B_2$ is defined to be the set
$B_1\otimes B_2=\{\,b_1\otimes b_2\,|\,b_i\in B_i\,\, (i=1,2)\,\}$
with
\begin{equation*}
{\rm wt}(b_1\otimes b_2)={\rm wt}(b_1)+{\rm wt}(b_2),
\end{equation*}
\begin{equation*}
\varepsilon_i(b_1\otimes b_2)= {\rm
max}(\varepsilon_i(b_1),\varepsilon_i(b_2)-\langle {\rm
wt}(b_1),h_i\rangle),
\end{equation*}
\begin{equation*}
\varphi_i(b_1\otimes b_2)= {\rm max}(\varphi_i(b_1)+\langle {\rm
wt}(b_2),h_i\rangle,\varphi_i(b_2)),
\end{equation*}
\begin{equation*}
{\te}_i(b_1\otimes b_2)=
\begin{cases}
{\te}_i b_1 \otimes b_2, & \text{if $\varphi_i(b_1)\geq \varepsilon_i(b_2)$}, \\
b_1\otimes {\te}_i b_2, & \text{if
$\varphi_i(b_1)<\varepsilon_i(b_2)$},
\end{cases}
\end{equation*}
\begin{equation*}
{\tf}_i(b_1\otimes b_2)=
\begin{cases}
{\tf}_i b_1 \otimes b_2, & \text{if  $\varphi_i(b_1)>\varepsilon_i(b_2)$}, \\
b_1\otimes {\tf}_i b_2, & \text{if $\varphi_i(b_1)\leq
\varepsilon_i(b_2)$},
\end{cases}
\end{equation*}
for $i\in S^\circ$ and $b_1\otimes b_2\in B_1\otimes B_2$, where we
assume that ${\bf 0}\otimes b_2=b_1\otimes {\bf 0}={\bf 0}$. Then
$B_1\otimes B_2$ is a crystal. If $B_1$ and $B_2$ are regular, then
so is $B_1\otimes B_2$. Note that $(B_1\otimes B_2)^\vee\simeq
B_2^\vee\otimes B_1^\vee$.

For a crystal $B$ and $m\in\Z_{\geq 0}$, we denote by $B^{\oplus m}$
the disjoint union $B_1\sqcup\cdots\sqcup B_m$ with $B_i\simeq B$,
where $B^{\oplus 0}$ means the empty set.

\subsection{}
Let $U_q(\gl_S)$ be the quantized enveloping algebra of $\gl_S$. For
$\Lambda\in P$, let $\B(\Lambda)$ be the crystal base of the
extremal weight $U_q(\gl_{S})$-module with extremal weight vector
$u_{\Lambda}$ of weight $\Lambda$, which is a regular
$\gl_S$-crystal.  We refer the reader to \cite{Kas94',Kas02} for
more details. When $\pm\Lambda$  is a dominant weight for $\gl_{S}$,
$\B(\Lambda)$ is a crystal base of the integrable highest (resp.
lowest) weight $U_q(\gl_{S})$-module with highest (resp. lowest)
weight $\Lambda$. Also we have $\B(\Lambda)\simeq \B(w\Lambda)$ for
$w\in W$. Hence, when ${S}$ is finite,  $\B(\Lambda)$ is always
isomorphic to the crystal base of a highest weight module and in
particular it is connected. When ${S}$ is infinite, it is shown in
\cite[Proposition 3.1]{K09} and \cite[Proposition 4.1]{K09-2} that
$\B(\Lambda)$ is also connected.

Let $\B(\pm \infty)$ be the crystal base of the negative (resp.
positive) part of $U_q(\gl_S)$ with the highest (resp. lowest)
weight vector $u_{\pm\infty}$ and let $T_\Lambda=\{\,t_\Lambda\,\}$
($\Lambda\in P$) be the crystal with $\te_i t_\Lambda=\tf_i
t_\Lambda ={\bf 0}$, ${\rm wt}(t_\Lambda)=\Lambda$ and
$\varepsilon_i(t_\Lambda)=\varphi_i(t_\Lambda)=-\infty$ for $i\in
S^\circ$. There is a strict embedding of $\B(\Lambda)$ into
$\B(\infty)\otimes T_\Lambda\otimes \B(-\infty)$ sending
$u_\Lambda$ to $u_\infty\otimes t_\Lambda\otimes u_{-\infty}$.
Hence $\B(\Lambda)$ is isomorphic to $C(u_\infty\otimes
t_\Lambda\otimes u_{-\infty})$  since $\B(\Lambda)$ is connected.
Note that $\B(\infty)\otimes T_\Lambda\otimes \B(-\infty)$ is regular.

There is an embedding  $\B(\Lambda_+)
\rightarrow \B(\infty)\otimes T_{\Lambda_+}$ (resp.
$\B(-\Lambda_-) \rightarrow T_{\Lambda_-}\otimes \B(-\infty)$)
sending $u_{\lambda_+}$ to $u_{\infty}\otimes t_{\Lambda_+}$ (resp.
$u_{\lambda_-}$ to $t_{\Lambda_-}\otimes u_{-\infty}$). This gives a
strict embedding
\begin{equation}
\iota_{\Lambda_+,\Lambda_-} : \B(\Lambda_+)\otimes \B(-\Lambda_-) \longrightarrow \B(\infty)\otimes
T_\Lambda\otimes \B(-\infty)
\end{equation}
sending $u_{\Lambda_+}\otimes u_{-\Lambda_-}$ to $u_\infty\otimes
t_\Lambda\otimes u_{-\infty}$ since $t_\Lambda\equiv
t_{\Lambda_+}\otimes t_{\Lambda_-}$. For a $\gl_S$-dominant weight
$\xi \in P$, let
\begin{equation}\label{embedding}
\iota_{\Lambda_+,\Lambda_-}^\xi :
\B(\Lambda_+)\otimes \B(-\Lambda_-) \longrightarrow \B(\Lambda_++\xi)\otimes \B(-\xi-\Lambda_-)
\end{equation}
be an embedding given by the composition of the following two morphisms
\begin{equation*}
\begin{split}
\B(\Lambda_+)\otimes \B(-\Lambda_-) & \longrightarrow
\B(\Lambda_+)\otimes \B(\xi)\otimes \B(-\xi)\otimes  \B(-\Lambda_-) \\
& \longrightarrow \B(\Lambda_++\xi)\otimes \B(-\xi-\Lambda_-),
\end{split}
\end{equation*}
where
\begin{equation*}
\begin{split}
&\tf_{i_1}\cdots\tf_{i_r}u_{\Lambda_+}\otimes
\te_{j_1}\cdots\te_{j_s}u_{-\Lambda_-} \\
&\mapsto \left(\tf_{i_1}\cdots\tf_{i_r}u_{\Lambda_+}\right)\otimes
u_{\xi}\otimes u_{-\xi}\otimes
\left(\te_{j_1}\cdots\te_{j_s}u_{-\Lambda_-}\right)\\
&\mapsto \tf_{i_1}\cdots\tf_{i_r}u_{\Lambda_++\xi}\otimes
\te_{j_1}\cdots\te_{j_s}u_{-\xi-\Lambda_-}.
\end{split}
\end{equation*}
Note that $$\tf_{i_1}\cdots\tf_{i_r}u_{\Lambda_++\xi}\equiv
\left(\tf_{i_1}\cdots\tf_{i_r}u_{\Lambda_+}\right)\otimes u_{\xi},\ \
\text{if $\tf_{i_1}\cdots\tf_{i_r}u_{\Lambda_+}\neq {\bf 0}$},$$
$$\te_{j_1}\cdots\te_{j_s}u_{-\xi-\Lambda_-}\equiv
u_{-\xi}\otimes \left(\te_{j_1}\cdots\te_{j_s}u_{-\Lambda_-}\right),\
\ \text{if $\te_{j_1}\cdots\te_{j_s}u_{-\Lambda_-}\neq {\bf 0}$.}$$
Since
\begin{equation}\label{characterization of Btilde}
\begin{split}
\B(\infty)\otimes T_{\Lambda}\otimes \B(-\infty)
&=\bigcup_{\substack{\text{$\Lambda', \Lambda''$ : $\gl_S$-dominant} \\
\Lambda'-\Lambda''=\Lambda}} {\rm Im}(\iota_{\Lambda',\Lambda''}),\\
\iota_{\Lambda',\Lambda''}&=\iota_{\Lambda'+\xi,\Lambda''+\xi}\circ
\iota^\xi_{\Lambda',\Lambda''},
\end{split}
\end{equation}
$\{\,\B(\Lambda')\otimes \B(-\Lambda'')\,|\,\text{$\Lambda',
\Lambda''$ : $\gl_S$-dominant with
$\Lambda=\Lambda'-\Lambda''$}\,\}$ forms a direct system, whose
limit is $\B(\infty)\otimes T_\Lambda\otimes \B(-\infty)$. Note
that $\B(\Lambda)$  is isomorphic to $C(u_{\Lambda_++\xi}\otimes
u_{-\xi-\Lambda_-})$ in
$\B(\Lambda_++\xi)\otimes \B(-\xi-\Lambda_-)$ for any $\gl_S$-dominant weight $\xi$.\vskip
2mm

\section{Young and Littlewood-Richardson tableaux}
\subsection{}
Let $\cP$ denote the set of partitions. We identify a partition
$\lambda=(\lambda_i)_{i\geq 1}$ with a {\it Young diagram} or a
subset $\{\,(i,j)\,|\,1\leq j\leq \lambda_i\,\}$ of
$\mathbb{N}\times\mathbb{N}$ following \cite{Mac95}. Let
$\ell(\lambda)=\left|\{\,i\,|\,\lambda_i\neq 0\,\}\right|$ and
$|\lambda|=\sum_{i\geq 1}\lambda_i$. We denote by
$\lambda'=(\lambda'_i)_{i\geq 1}$ the conjugate partition of
$\lambda$ whose Young diagram is
$\{\,(i,j)\,|\,(j,i)\in\lambda\,\}$. For $\mu,\nu\in\cP$,
$\mu\cup\nu$ is the partition obtained  by rearranging $\{\,\mu_i,
\nu_i\,|\,i\geq 1\,\}$, and $\mu+\nu=(\mu_i+\nu_i)_{i\geq 1}$.

Let $\A$ be a linearly ordered set and $\lambda/\mu$ a skew Young
diagram. A tableau $T$ obtained by filling $\lambda/\mu$ with
entries in $\A$ is called a {\it semistandard tableau or Young
tableau of shape $\lambda/\mu$} if the entries in each row (resp.
column) are weakly (resp. strictly) increasing from left to right
(resp. from top to bottom). We denote by $T(i,j)$ the entry of $T$
at $(i,j)$ position for $(i,j)\in \lambda/\mu$. Let
$SST_\A(\lambda/\mu)$ be the set of all semistandard tableaux of
shape $\lambda/\mu$ with entries in $\A$.

Suppose that $\A$ is an interval in $\Z$ with a usual linear
ordering. Then $\A$ is a $\gl_\A$-crystal associated with the
natural representation of $U_q(\gl_\A)$, where ${\rm
wt}(i)=\epsilon_i$  and $i \stackrel{i}{\rightarrow} i+1$ for $i\in
\A^\circ$. The image of $SST_\A(\lambda/\mu)$ in $\A^{\otimes r}$
under the map $T \mapsto w(T)_{\rm col}=w_1\ldots w_r$ together with
$\{{\bf 0}\}$ is invariant under $\te_i,\tf_i$, where $w(T)_{\rm
col}$  is the word obtained by reading the entries column by column
from right to left, and in each column from top to bottom. Hence
$SST_\A(\lambda/\mu)$ is a subcrystal of $\A^{\otimes r}$ \cite{KN}.
We may identify $T^\vee \in SST_\A(\lambda/\mu)^\vee$ with the
tableau obtained from $T$ by $180^\circ$-rotation and replacing each
entry $t$ with $t^\vee$. So we have $SST_\A(\lambda/\mu)^\vee\simeq
SST_{\A^\vee}((\lambda/\mu)^\vee)$, where $a^\vee < b^\vee$ if and
only if $b<a$ for $a,b\in\A$ and $(\lambda/\mu)^\vee$ is the skew
Young diagram obtained from $\lambda/\mu$ by $180^\circ$-rotation.
We use the convention $(t^{\vee})^{\vee}=t$ and hence
$(T^{\vee})^{\vee}=T$.

\subsection{}
For $\lambda,\mu,\nu\in\cP$ with $|\lambda|=|\mu|+|\nu|$, let ${\rm\bf LR}^\lambda_{\mu
\nu}$ be the set of tableaux $U$ in
$SST_{\mathbb{N}}(\lambda/\mu)$ such that
\begin{itemize}
\item[(1)] the number of occurrences of each $i\geq 1$ in $U$ is
$\nu_i$,

\item[(2)] for for $1\leq k\leq |\nu|$, the number of occurrences of each $i\geq 1$ in
$w_1\ldots w_k$ is no less than that of $i+1$ in $w_1\ldots w_k$,
where $w(U)_{\rm col}=w_1\ldots w_{|\nu|}$.
\end{itemize}
We call ${\rm\bf LR}^\lambda_{\mu \nu}$ the set of {\it
Littlewood-Richardson tableaux of shape $\lambda/\mu$ with content
$\nu$} and put $c^\lambda_{\mu \nu}=\left|{\rm\bf LR}^\lambda_{\mu
\nu}\right|$ \cite{Mac95}. We introduce a variation of ${\rm\bf
LR}^\lambda_{\mu \nu}$,  which is necessary for our later arguments.
Let $\overline{{\rm\bf LR}}^\lambda_{\mu \nu}$ be the set of
tableaux $U$ in $SST_{-\mathbb{N}}(\lambda/\mu)$ such that
\begin{itemize}
\item[(1)] the number of occurrences of each $-i\leq -1$ in $U$ is
$\nu_i$,

\item[(2)] for $1\leq k\leq |\nu|$, the number of occurrences of each $-i\leq -1$ in
$w_k\ldots w_{|\nu|}$ is no less than that of $-(i+1)$ in $w_k\ldots
w_{|\nu|}$, where $w(U)_{\rm col}=w_1\ldots w_{|\nu|}$.
\end{itemize}

There are other characterizations of ${\rm\bf LR}^\lambda_{\mu \nu}$
and $\overline{{\rm\bf LR}}^\lambda_{\mu \nu}$ using crystals. For
$U\in SST_{\mathbb{N}}(\lambda/\mu)$, we have $U\in {\rm\bf
LR}^\lambda_{\mu \nu}$ if and only if $U$ is $\gl_{>0}$-equivalent
(or Knuth equivalent) to the highest weight element $H_{\nu}$ in
$SST_{\mathbb{N}}(\nu)$, that is, $H_\nu(i,j)=i$ for $(i,j)\in \nu$.
Similarly, we have for $U\in SST_{-\mathbb{N}}(\lambda/\mu)$, $U\in
\overline{{\rm\bf LR}}^\lambda_{\mu \nu}$ if and only if $U$ is
$\gl_{<0}$-equivalent (or Knuth equivalent) to the lowest weight
element $L_{\nu}$ in $SST_{-\mathbb{N}}(\nu)$, that is,
$L_\nu(i,j)=-\nu'_j+i-1$ for $(i,j)\in\nu$.

There is a one-to-one correspondence with the set of $V\in
SST_{\mathbb{N}}(\nu)$ such that $H_\mu \otimes V \equiv H_\lambda$
and ${\rm\bf LR}^\lambda_{\mu \nu}$. Indeed, $V$ corresponds to
$\imath(V)=U\in {\rm\bf LR}^\lambda_{\mu \nu}$ where the number of
$k$'s in the $i$-th row of $V$ is equal to the number of $i$'s in
the $k$-th row of $U$ for $i, k\geq 1$.

\begin{ex}{\rm\mbox{}
\begin{center}
$SST_{\mathbb{N}}((3,3,2)) \ni$
\begin{tabular}{ccc}
    1 & 1 & 2  \\
    2 & 2 & 3  \\
    3 & 4
\end{tabular} \ \
\ \ $\stackrel{\imath}{\longrightarrow}$ \ \ \ \
\begin{tabular}{ccccc}
    {$\bullet$} & {$\bullet$} & {$\bullet$} &  1 & 1 \\
    {$\bullet$} & 1 & 2 & 2 & \\
    2 & 3 & & & \\
    3 &  & & &
\end{tabular}
$\in {\bf LR}^{(5,4,2,1)}_{(3,1)\ (3,3,2)}$
\end{center}\vskip 2mm}
\end{ex}

\subsection{}
Next, let us briefly recall the {\it switching algorithm} \cite{BSS}. Suppose that $\mathcal{A}$ and
$\mathcal{B}$ are two linearly ordered sets. Let $\lambda/\mu$ be a skew Young diagram.
Let $U$ be a tableau of shape $\lambda/\mu$ with
entries in $\mathcal{A}\sqcup\mathcal{B}$, satisfying the following
conditions;
\begin{itemize}
\item[(S1)] $U(i,j)\leq U(i',j')$ whenever $U(i,j), U(i',j')\in\mathcal{X}$
for $(i,j), (i',j')\in\lambda/\mu$ with $i\leq i'$ and $j\leq j'$,

\item[(S2)] in each column of $U$,
entries in $\mathcal{X}$  increase strictly from top to bottom,
\end{itemize}
where  $\mathcal{X}=\A$ or $\mathcal{B}$. Suppose that $a\in \mathcal{B}$ and $b\in
\mathcal{A}$ are two adjacent entries in $U$ such that $a$ is placed
above or to the left of $b$. Interchanging $a$ and $b$ is called a
{\it switching} if the resulting tableau still satisfies the
conditions (S1) and (S2).

Let $\lambda/\mu$ and $\mu/\eta$ be two skew Young diagrams.  For
$S\in SST_{\mathcal{B}}(\mu/\eta)$ and $T\in
SST_{\mathcal{A}}(\lambda/\mu)$, we denote by $S\ast T$ the tableau
of shape $\lambda/\eta$ with entries
${\mathcal{A}\sqcup\mathcal{B}}$ obtained by gluing $S$ and $T$. Let
$U$ be a tableau obtained from $S\ast T$ by applying switching
procedures as far as possible. Then it is shown in \cite[Theorems
2.2 and 3.1]{BSS} that
\begin{itemize}
\item[(1)] $U=T'\ast S'$,
where $T'\in SST_{\A}(\nu/\eta)$ and $S'\in
SST_{\mathcal{B}}(\lambda/\nu)$ for some $\nu$,

\item[(2)] $U$ is uniquely determined by $S$ and $T$,

\item[(3)] $w(S)_{\rm col}$ (resp. $w(T)_{\rm col}$) is Knuth equivalent
to $w(S')_{\rm col}$ (resp. $w(T')_{\rm col}$),

\end{itemize}
Suppose that $\eta=\emptyset$ and $S=H_{\mu}\in
SST_{\mathbb{N}}(\mu)$. We put
\begin{equation}
\jmath(T)=T', \ \ \ \ \ \jmath(T)_R=S'.
\end{equation}
Then we have the following.
\begin{prop}\label{skew LR}
Suppose that $\A$ is an interval in $\Z$. The map sending $T$ to
$(\jmath(T),\jmath(T)_R)$ is an isomorphism of $\gl_\A$-crystals
$$SST_{\A}(\lambda/\mu) \longrightarrow \bigsqcup_{\nu\in\cP}
SST_{\A}(\nu)\times {\bf LR}^{\lambda}_{\nu\,\mu},$$ where
$\widetilde{x}_i(T',S')=(\widetilde{x}_iT',S')$ for $i\in \A^\circ$
and $x=e,f$ on the righthand side. In particular, the map $Q \mapsto
\jmath(Q)_R$ restricts to a bijection from ${\bf
LR}^{\lambda}_{\mu\,\nu}$ to ${\bf LR}^{\lambda}_{\nu\,\mu}$, and
from $\overline{\bf LR}^{\lambda}_{\mu\,\nu}$ to ${\bf
LR}^{\lambda}_{\nu\,\mu}$ when $\A=\pm\mathbb{N}$, respectively.
\end{prop}
\pf The map is clearly a bijection by \cite[Theorem 3.1]{BSS}.
Moreover, $\jmath(T)$ is $\gl_{\A}$-equivalent to $T$ and
$\jmath(T)_R$ is invariant under $\te_i$ and $\tf_i$ for $i\in
\A^\circ$ (cf.\cite[Theorem 5.9]{HK}). Hence the bijection is an
isomorphism of $\gl_\A$-crystal, where on the right-hand side the
operators $\te_i$, $\tf_i$ act on the first component. \qed

\begin{rem}{\rm
The inverse of the isomorphism in Proposition \ref{skew LR}
is given directly by applying the switching process in
a reverse way.}
\end{rem}

\section{Extremal weight crystals of type $A_{+\infty}$}
Note that for $r\in\Z$ the $\gl_{>r}$-crystals $[r+1,\infty)$ and $[r+1,\infty)^\vee$ are given by
{\allowdisplaybreaks
\begin{equation*}\label{natural representations}
\begin{split}
 \ \ &\ \ \  r+1 \stackrel{\tiny r+1}{\longrightarrow} r+2
\stackrel{r+2}{\longrightarrow}r+3\stackrel{r+3}{\longrightarrow}
\cdots,\\
\ \ & \cdots \stackrel{r+3}{\longrightarrow}
(r+3)^\vee\stackrel{r+2}{\longrightarrow}(r+2)^\vee\stackrel{r+1}{\longrightarrow}
(r+1)^\vee  \ .
\end{split}
\end{equation*}}
For $\mu\in \cP$, let
\begin{equation}
\begin{split}
\B_{\mu}^{>r}= SST_{[r+1,\infty)}(\mu). \\
\end{split}
\end{equation}
Then $\B_{\mu}^{>r}$ is the highest weight $\gl_{>r}$-crystal with
highest weight element  $H^{>r}_\mu$ of weight $\sum_{i\geq
1}\lambda_i\epsilon_{r+i}$, where $H^{>r}_\mu(i,j)=r+i$ for
$(i,j)\in\mu$. We identify $\left(\B_{\mu}^{>r}\right)^\vee$ with
$SST_{[r+1,\infty)^\vee}(\mu^\vee)$.

For $\nu\in \cP$ and $s\geq \ell(\nu)$, let $E^{>r}_{\nu}(s)\in
\left(\B_{\nu}^{>r}\right)^\vee$ be an element given by
\begin{equation}\label{Extremal weight vector}
\left(E^{>r}_{\nu}(s)\right)^\vee(i,j)=s-\nu_j'+i
\end{equation}
for $(i,j)\in\nu$. For $s\geq \ell(\mu)+\ell(\nu)$, let
\begin{equation}
\B_{\mu,\nu}^{> r} = C\left(H^{>r}_\mu\otimes E^{>r}_{\nu}(s)\right)
\subset \B_{\mu}^{>r}\otimes \left(\B_{\nu}^{>r}\right)^\vee,
\end{equation}
the connected component including $H^{>r}_\mu\otimes
E^{>r}_{\nu}(s)$ as a $\gl_{>r}$-crystal.  Then we have the
following by \cite[Proposition 3.4]{K09} and \cite[Theorem
3.5]{K09}.

\begin{thm}\label{extremal +infty} For $\mu,\nu\in\cP$,
\begin{itemize}
\item[(1)] $\B_{\mu,\nu}^{>r}$ is the set of $S\otimes
T\in \B_{\mu}^{>r}\otimes \left(\B_{\nu}^{>r}\right)^\vee$ such that
for each $k\geq 1$,
$$\Bigl|\,\{\,i\,|\,S(i,1)\leq r+k\,\}\, \Bigr|+\left|\,\{\,i\,|\,T^\vee(i,1)\leq r+k\,\}\, \right| \leq k,$$

\item[(2)] $\B_{\mu,\nu}^{> r}$ is isomorphic to
an extremal weight $\gl_{>r}$-crystal with extremal weight
$$\sum_{i=1}^{\ell(\mu)}\mu_i\epsilon_{r+i}-\sum_{j=1}^{\ell(\nu)}\nu_j\epsilon_{\ell(\mu)+\ell(\nu)-j+1}.$$
\end{itemize}
\end{thm}

Note that $\B_{\mu,\nu}^{> r}$ does not depend on the choice of $s$
and $\{\,\B_{\mu,\nu}^{> r}\,|\,\mu,\nu\in\cP\,\}$  is a complete
list of pairwise non-isomorphic extremal weight $\gl_{>r}$-crystals
\cite[Theorem 3.5 and Lemma 5.1]{K09} and the tensor product of
extremal weight $\gl_{>r}$-crystals is isomorphic to a finite
disjoint union of extremal weight crystals \cite[Theorem 4.10]{K09}.

To describe the tensor product of extremal weight $\gl_{>r}$-crystals,
let us review an insertion algorithm for extremal weight crystal
elements \cite{K09}, which is an infinite analogue of \cite{Str}.
Recall that for $a\in\A$ and $T\in SST_\A(\lambda)$, $a\rightarrow
T$ (resp. $T\leftarrow a$) denotes the tableau obtained by the
Schensted column (resp. row) insertion, where $\A$ is a linearly
ordered set and $\lambda\in\cP$ (see for example \cite[Appendix
A.2]{Fu}).

We denote $S\otimes T\in\B^{>r}_{\mu,\nu}$ by $(S,T)$. For $a\in
[r+1,\infty)$, we define $a\rightarrow (S,T)$ in the following way;

Suppose first that $S$ is empty and $T$ is a single column tableau. Let
$(T',a')$ be the pair obtained by the following process;

\begin{itemize}
\item[(1)] If $T$ contains $a^\vee, (a+1)^\vee,\ldots, (b-1)^\vee$
but not $b^\vee$, then $T'$ is the tableau obtained from $T$ by
replacing $a^\vee, (a+1)^\vee,\ldots, (b-1)^\vee$ with
$(a+1)^\vee,(a+2)^\vee,\ldots, b^\vee$, and put $a'=b$.

\item[(2)] If $T$ does not contain $a^\vee$, then leave $T$
unchanged and put $a'=a$.
\end{itemize}
Now, we suppose that $S$ and $T$ are arbitrary.
\begin{itemize}
\item[(1)] Apply the above process to the leftmost column of $T$ with $a$.

\item[(2)] Repeat (1) with $a'$  and the next column to the right.

\item[(3)] Continue this process to the right-most column of $T$ to get a tableau $T'$ and $a''$.

\item[(4)] Define
$a\rightarrow (S,T)$ to be $\left((a''\rightarrow S)\,,\,T'\right)$
\end{itemize}
Then $(a\rightarrow (S,T))\in\B^{>r}_{\sigma,\nu}$ for some $\sigma
\in\cP$ with $|\sigma|-|\mu|=1$. For a finite word $w=w_1\ldots w_n$
with letters in $[r+1,\infty)$, we let $(w\rightarrow (S,T))=(
w_n\rightarrow(\cdots(w_1\rightarrow (S,T))\cdots))$.

For $a\in [r+1,\infty)$ and $(S,T)\in\B^{>r} _{\mu,\nu}$, we define
$(S,T)\leftarrow a^\vee$ to be the pair $(S',T')$ obtained in the
following way;
\begin{itemize}
\item[(1)] If the pair $(S,(T^\vee \leftarrow a)^\vee)$
satisfies the condition in Theorem \ref{extremal +infty} (1),
then put $S'=S$ and $T'=(T^\vee \leftarrow a)^\vee$.

\item[(2)] Otherwise, choose the smallest $k$
such that $a_k$ is bumped out of the $k$-th row in the row insertion
of $a$ into $T^\vee$ and the insertion of $a_k$ into the $(k+1)$-th
row violates the condition in Theorem \ref{extremal +infty} (1).

\item[(2-a)] Stop the row insertion of $a$ into
$T^\vee$ when $a_k$ is bumped out and let $T'$ be the resulting
tableau after taking $\vee$.

\item[(2-b)] Remove $a_k$ in the left-most column of $S$, which
necessarily exists, and then apply the {\it jeu de taquin} (see for
example \cite[Section 1.2]{Fu}) to obtain a tableau $S'$.
\end{itemize}
In this case, $((S,T)\leftarrow a^\vee)\in \B^{>r}_{\sigma,\tau}$,
where either (1) $|\mu|-|\sigma|=1$ and $\tau=\nu$, or (2)
$\sigma=\mu$ and $|\tau|-|\nu|=1$. For a finite word $w=w_1\ldots
w_n$ with letters in $[r+1,\infty)^\vee$, we let $((S,T)\leftarrow
w)=((\cdots((S,T)\leftarrow w_1)\cdots )\leftarrow w_n)$.

Let $\mu,\nu,\sigma,\tau\in\cP$ be given. For
$(S,T)\in\B^{>r}_{\mu,\nu}$ and $(S',T')\in\B^{>r}_{\sigma,\tau}$,
we define
\begin{equation*}
\left((S',T')\rightarrow (S,T)\right) = \left(\left( w(S')_{\rm
col}\rightarrow (S,T)\right) \leftarrow w(T')_{\rm col}\right).
\end{equation*}
Then $((S',T')\rightarrow (S,T))\in\B^{>r}_{\zeta,\eta}$ for some
$\zeta,\eta\in\cP$.   Assume that $w(S')_{\rm col}=w_1\ldots
w_s$ and $w(T')_{\rm col}=w_{s+1}\ldots w_{s+t}$. For $1\leq i\leq
s+t$, let
\begin{equation*}
(S^i,T^i)=
\begin{cases}
w_1\cdots w_i\rightarrow (S,T), & \text{if $1\leq i\leq
s$}, \\
(S^s,T^s)\leftarrow w_{s+1}\cdots w_i, & \text{if $s+1\leq i\leq
s+t$},
\end{cases}
\end{equation*}
and $(S^0,T^0)=(S,T)$. We define
\begin{equation*}
\left((S',T')\rightarrow (S,T)\right)_R=(U,V),
\end{equation*}
where $(U,V)$ is the pair of tableaux with entries in
$\mathbb{Z}\setminus\{0\}$ determined by the following process;
\begin{itemize}
\item[(1)] $U$ is of shape $\sigma$ and $V$ is of shape $\tau$.

\item[(2)] Let $1\leq i\leq s$. If $w_i$ is inserted into
$(S^{i-1},T^{i-1})$ to create a dot (or box) in the $k$-th row of
the shape of $S^{i-1}$, then we fill the dot in $\sigma$
corresponding to $w_i$ with $k$.

\item[(3)] Let $s+1\leq i\leq s+t$. If $w_i$ is inserted into
$(S^{i-1},T^{i-1})$ to create a dot in the $k$-th row (from the
bottom) of the shape of $T^{i-1}$, then we fill the dot in $\tau$
corresponding to $w_i$ with $-k$. If $w_i$ is inserted into
$(S^{i-1},T^{i-1})$ to remove a dot in the $k$-th row of the shape
of $S^{i-1}$, then we fill the corresponding dot in $\tau$ with $k$.
\end{itemize}
We call $\left((S',T')\rightarrow (S,T)\right)_R$ {\it the recording
tableau of $\left((S',T')\rightarrow (S,T)\right)$}. By
\cite[Theorem 4.10]{K09}, we have the following.
\begin{prop}\label{recording} Under the above hypothesis, we have
\begin{itemize}
\item[(1)] $\left((S',T')\rightarrow (S,T)\right)\equiv (S,T)\otimes
(S',T')$,

\item[(2)] $\left((S',T')\rightarrow (S,T)\right)_R\in SST_\mathbb{N}(\sigma)\times
SST_{\mathcal{Z}}(\tau)$, where  $\mathcal{Z}$ is the set of non-zero integers with a linear
ordering $1\prec2\prec3\prec\cdots\prec-3\prec-2\prec-1$,

\item[(3)] the recording tableaux are constant on the connected component of $\B^{>r}_{\mu,\nu}\otimes
\B^{>r}_{\sigma,\tau}$  including $(S,T)\otimes (S',T')$.
\end{itemize}
\end{prop}

Suppose that $\mu,\nu\in\cP$ and $W\in SST_\mathcal{Z}(\nu)$ are
given with $w(W)_{\rm col}=w_n\ldots w_1$. Let
$(\alpha^0,\beta^0),(\alpha^1,\beta^1),\ldots,(\alpha^{n},\beta^{n})$
be the sequence where $\alpha^i=(\alpha^i_j)_{j\geq 1}$ and
$\beta^i=(\beta^i_j)_{j\geq 1}$ ($1\leq i\leq n$) are sequences of
integers defined inductively as follows;
\begin{itemize}
\item[(1)] $\alpha^0=\mu$ and $\beta^0=\emptyset$.

\item[(2)] If $w_i$ is positive, then
$\alpha^{i}$ is obtained by subtracting $1$ in the $w_i$-th part of
$\alpha^{i-1}$, and $\beta^{i}=\beta^{i-1}$. If $w_i$ is negative,
then $\alpha^{i}=\alpha^{i-1}$ and $\beta^{i}$ is obtained by adding
$1$ in the $(-w_i)$-th part of $\beta^{i-1}$.
\end{itemize}
Then for $\sigma, \tau\in\cP$ we define
$\mathcal{C}^{(\mu,\nu)}_{(\sigma,\tau)}$ to be the set of $W\in
SST_\mathcal{Z}(\nu)$ such that (1) $\alpha^i$, $\beta^i\in \cP$ for
$1\leq i\leq n$, and (2) $\alpha^{n}=\sigma$, $\beta^{n}=\tau$.

As a particular case of \cite[Theorem 4.10]{K09}, we have the following.

\begin{prop} For $\mu,\nu\in\cP$, we have an isomorphism of $\gl_{>r}$-crystals
\begin{equation*}
\B_{\mu}^{>r}\otimes \left(\B_{\nu}^{>r}\right)^\vee \longrightarrow
\bigsqcup_{\sigma,\tau\in\cP}\B_{\sigma,\tau}^{>r}\times
\mathcal{C}^{(\mu,\nu)}_{(\sigma,\tau)},
\end{equation*}
where $S\otimes T$ is sent to $(\left((\emptyset,T)\rightarrow
(S,\emptyset)\right),\left((\emptyset,T)\rightarrow
(S,\emptyset)\right)_R)$.
\end{prop}

Further, we can characterize
$\mathcal{C}^{(\mu,\nu)}_{(\sigma,\tau)}$ as follows.
\begin{prop}\label{LR decomp} For $\mu,\nu,\sigma,\tau\in\cP$,  there exists a bijection
\begin{equation*}
\mathcal{C}^{(\mu,\nu)}_{(\sigma,\tau)} \longrightarrow
\bigsqcup_{\lambda\in\cP}{\bf LR}^\mu_{\sigma \lambda}\times {\bf
LR}^{\nu}_{\tau\lambda }.
\end{equation*}
\end{prop}
\pf Suppose that $W\in \mathcal{C}^{(\mu,\nu)}_{(\sigma,\tau)}$ is
given. Let $W_{+}$ (resp. $W_{-}$) be the subtableau in $W$
consisting of positive (resp. negative) entries.

We have  $W_{+}\in SST_{\mathbb{N}}(\lambda)$ and $W_{-}\in
SST_{-\mathbb{N}}(\nu/\lambda)$ for some $\lambda\subset \nu$. By
definition of $W\in \mathcal{C}^{(\mu,\nu)}_{(\sigma,\tau)}$,  we
have $\imath(W_{+})\in {\bf LR}^{\mu}_{\sigma \lambda}$ and
$W_{-}\in \overline{\bf LR}^{\nu}_{\lambda \tau}$, hence
$\jmath(W_{-})_R\in {\bf LR}^{\nu}_{\tau\lambda }$ by Proposition
\ref{skew LR}.

We can check that the correspondence
$W\mapsto(W_1,W_2)=(\imath(W_{+}),\jmath(W_{-})_R)$ is reversible
and hence gives a bijection $\mathcal{C}^{(\mu,\nu)}_{(\sigma,\tau)}
\longrightarrow \bigsqcup_{\lambda\in\cP}{\bf LR}^\mu_{\sigma
\lambda}\times {\bf LR}^{\nu}_{\tau\lambda }$. \qed\vskip 2mm

\begin{ex}\label{insertion ex-1}{\rm Consider
$$
S=
\begin{array}{ccc}
1 & 1   & 2  \\
2 & 3
\end{array} \in \B^{>0}_{(3,2)}, \ \ \ \
T=
\begin{array}{ccc}
    &  &  4^\vee  \\
   & 3^\vee &  2^\vee  \\
 2^\vee   & 2^\vee  & 1^\vee
\end{array} \in \left(\B^{>0}_{(3,2,1)}\right)^\vee.
$$
Then we have
$$\left(
\begin{array}{ccc}
1 & 1   & 2  \\
2 & 3 &
\end{array}, \ \
\ \ \emptyset \ \  \right)\ \longleftarrow\ 4^\vee =\
\left(\begin{array}{ccc}
1 & 1   & 2  \\
2 & 3
\end{array}, \ \
\begin{array}{cc}
     &    \\
     &  4^\vee  \\
 \end{array}
\right)\ \ \ \ \ \ \ \ \ \ \
\begin{array}{ccc}
\bullet & \bullet   & \bullet  \\
\bullet & \bullet \\
\!\!\!-1
\end{array}$$
$$\left(\begin{array}{ccc}
1 & 1   & 2  \\
2 & 3 &
\end{array}, \ \
\begin{array}{cc}
     &    \\
     &  4^\vee  \\
 \end{array}
\right)\ \longleftarrow\ 2^\vee =\ \left(
\begin{array}{ccc}
1 & 1   & 2  \\
3 &
\end{array}, \ \
\begin{array}{cc}
     &    \\
     &  4^\vee  \\
 \end{array}
\right)\ \ \ \ \ \ \ \ \ \ \
\begin{array}{ccc}
\bullet & \bullet   & \bullet  \\
2 & \bullet \\
\!\!\!-1
\end{array}$$
$$
\left(
\begin{array}{ccc}
1 & 1   & 2  \\
3 &
\end{array}, \ \
\begin{array}{cc}
     &    \\
     &  4^\vee  \\
 \end{array}
\right)\ \longleftarrow\ 1^\vee =\ \left(
\begin{array}{ccc}
1 & 2   &    \\
3 &
\end{array}, \ \
\begin{array}{cc}
     &    \\
     &  4^\vee  \\
 \end{array}
\right)\ \ \ \ \ \ \ \ \ \ \ \ \ \!
\begin{array}{ccc}
1 & \bullet   & \bullet  \\
2 & \bullet \\
\!\!\!-1
\end{array}$$
$$\left(
\begin{array}{ccc}
1 & 2   &    \\
3 &
\end{array}, \ \
\begin{array}{cc}
     &    \\
     &  4^\vee  \\
 \end{array}
\right)\ \longleftarrow\ 3^\vee =\ \left(
\begin{array}{ccc}
1 & 2   &    \\
3 &
\end{array}, \ \
\begin{array}{cc}
     &  4^\vee \\
     &  3^\vee  \\
 \end{array}
\right)\ \ \ \ \ \ \ \ \ \ \ \
\begin{array}{ccc}
1 & \bullet   & \bullet  \\
2 & \!\!\!-2 \\
\!\!\!-1
\end{array}$$
$$\left(
\begin{array}{ccc}
1 & 2   &    \\
3 &
\end{array}, \ \
\begin{array}{cc}
     &  4^\vee \\
     &  3^\vee  \\
 \end{array}
\right)\ \longleftarrow\ 2^\vee =\ \left(
\begin{array}{ccc}
1 & 2   &    \\
  &
\end{array}, \ \
\begin{array}{cc}
     &  4^\vee \\
     &  2^\vee  \\
 \end{array}
\right)\ \ \ \ \ \ \ \ \ \ \ \
\begin{array}{ccc}
1 & 2   & \bullet  \\
2 & \!\!\!-2 \\
\!\!\!-1
\end{array}$$
$$\left(
\begin{array}{ccc}
1 & 2   &    \\
  &
\end{array}, \ \
\begin{array}{cc}
     &  4^\vee \\
     &  2^\vee  \\
 \end{array}
\right) \ \longleftarrow\ 2^\vee =\ \left(
\begin{array}{ccc}
1 & 2   &    \\
  &
\end{array}, \ \
\begin{array}{cc}
     &  4^\vee \\
  2^\vee &  2^\vee  \\
 \end{array}
\right)\ \ \ \ \ \ \ \ \
\begin{array}{ccc}
1 & 2   & \!\!\!-1  \\
2 & \!\!\!-2 \\
\!\!\!-1
\end{array}.$$
Hence,
\begin{equation*}
\begin{split}
\left((\emptyset,T)\rightarrow (S,\emptyset)\right)&= \left(
\begin{array}{cc}
1 & 2       \\
  &
\end{array}, \ \
\begin{array}{cc}
     &  4^\vee \\
  2^\vee &  2^\vee  \\
 \end{array}
\right)\in \B^{>0}_{(2),(2,1)} ,\\
\left((\emptyset,T)\rightarrow (S,\emptyset)\right)_R &=
\begin{array}{ccc}
1 & 2   & \!\!\!-1  \\
2 & \!\!\!-2 \\
\!\!\!-1
\end{array}\in \mathcal{C}^{(3,2),(3,2,1)}_{(2), (2,1)}.
\end{split}
\end{equation*}
If we put $W=\left((\emptyset,T)\rightarrow (S,\emptyset)\right)_R$,
then
\begin{equation*}
W_{+}=
\begin{array}{ccc}
1 & 2   &   \\
2 &
\end{array}, \ \ \
W_{-}=
\begin{array}{ccc}
\bullet & \bullet   & \!\!\!-1  \\
\bullet & \!\!\!-2 \\
\!\!\!-1
\end{array}. \ \ \
\end{equation*}
Since
\begin{equation*}
\imath(W_{+}) =
\begin{array}{ccc}
\bullet & \bullet   & 1  \\
1 & 2
\end{array},\ \
{\jmath}(W_{-})=
\begin{array}{cc}
-2 & -1    \\
-1 &
\end{array}, \ \ \
{\jmath}(W_{-})_R=
\begin{array}{ccc}
\bullet & \bullet & 1   \\
\bullet & 2 \\
1
\end{array}
\end{equation*}
(see Proposition \ref{skew LR} (2)), we have
\begin{equation*}
(W_1, W_2)= \left(
\begin{array}{ccc}
\bullet & \bullet   & 1  \\
1 & 2 \\
&
\end{array}\ \ ,\ \
\begin{array}{ccc}
\bullet & \bullet & 1   \\
\bullet & 2 \\
1
\end{array}
\right)\in  {\bf LR}^{(3,2)}_{(2) (2,1)}\times {\bf
LR}^{(3,2,1)}_{(2,1) (2,1)}.
\end{equation*}
}
\end{ex}

\vskip 5mm

Now, the multiplicity of each connected component can be
written in terms of Littlewood-Richardson coefficients as follows. We remark that it was already given in
\cite[Corollary 7.3]{K09}, while Proposition \ref{LR decomp} gives a bijective proof of it.

\begin{cor}\label{LR rule for >r}
For $\mu,\nu\in\cP$, we have
\begin{equation*}
\B_{\mu}^{>r}\otimes \left(\B_{\nu}^{>r}\right)^\vee \simeq
\bigsqcup_{\sigma,\tau\in\cP}\left(\B_{\sigma,\tau}^{>r}\right)^{\oplus
c^{(\mu,\nu)}_{(\sigma,\tau)}},
\end{equation*}
where
\begin{equation*}
c^{(\mu,\nu)}_{(\sigma,\tau)}=\sum_{\lambda\in\cP}c^{\mu}_{
\sigma\lambda}c^{\nu}_{ \tau\lambda}.
\end{equation*}
\end{cor}

\begin{prop}\label{tensor decomposition} For $\mu,\nu\in\cP$, we have an isomorphism
of $\gl_{>r}$-crystals
$$\left(\B_{\nu}^{>r}\right)^\vee\otimes \B_{\mu}^{>r}
\longrightarrow \B^{>r}_{\mu,\nu},$$ where $T\otimes S$ is mapped to
$\left( (S,\emptyset)\rightarrow (\emptyset,T)\right)$.
\end{prop}
\pf For $T\otimes S\in \left(\B_{\nu}^{>r}\right)^\vee\otimes
\B_{\mu}^{>r}$, it follows from Proposition \ref{recording} (2) that
\begin{itemize}
\item[(1)] $\left( (S,\emptyset)\rightarrow
(\emptyset,T)\right)_R=\left(H_{\mu},\emptyset\right)$.

\item[(2)] $\left( (S,\emptyset)\rightarrow
(\emptyset,T)\right)\in \B^{>r}_{\mu,\nu}$,
\end{itemize}
Therefore, by  \cite[Theorem 4.10]{K09} the map
\begin{equation*}
\left(\B_{\nu}^{>r}\right)^\vee\otimes \B_{\mu}^{>r}
\stackrel{}{\longrightarrow} \B_{\mu,\nu}^{>r}\times
\{\,\left(H_\mu,\emptyset\right)\,\}
\end{equation*}
sending $T\otimes S$ to $\left(\left( (S,\emptyset)\rightarrow
(\emptyset,T)\right),\left( (S,\emptyset)\rightarrow
(\emptyset,T)\right)_R\right)$ is an isomorphism of
$\gl_{>r}$-crystals. \qed

\begin{ex}\label{insertion ex-2}{\rm Let
$$(U,V)= \left(
\begin{array}{cc}
1 & 2       \\
  &
\end{array}, \ \
\begin{array}{cc}
     &  4^\vee \\
  2^\vee &  2^\vee  \\
 \end{array}
\right)\in \B^{>0}_{(2),(2,1)}
$$
be as in  Example \ref{insertion ex-1}. If we put
$$
\widetilde{V}\otimes \widetilde{U} =
\begin{array}{cc}
     &  4^\vee \\
  2^\vee &  1^\vee  \\
 \end{array} \otimes
\begin{array}{cc}
1 & 1       \\
  &
\end{array} \ \
\in \left(\B_{(2,1)}^{>0}\right)^\vee\otimes \B_{(2)}^{>0},
$$
then $$\left((\widetilde{U},\emptyset)\rightarrow
(\emptyset,\widetilde{V})\right)=(U,V).$$ }
\end{ex}\vskip 5mm

\section{Combinatorial description of $\B(\widetilde{U}_q(\gl_{>0}))$}
In this section, we give a combinatorial realization of
$\B(\infty)\otimes T_\Lambda \otimes \B(-\infty)$ for an integral
weight $\Lambda$ in case of $\gl_{>0}$, and then
$\B(\widetilde{U}_q(\gl_{>0}))$.

\subsection{}
For simplicity, we put for a skew Young diagram $\lambda/\mu$
$$\cB_{\lambda/\mu}=SST_{\mathbb{N}}(\lambda/\mu),$$
and for $\mu,\nu\in\cP$ $$\cB_{\mu,\nu}=\B^{>0}_{\mu,\nu}.$$
For $S\otimes T\in \cB_\mu\otimes
\cB_\nu^\vee$, suppose that
\begin{equation*}
\begin{split}
(U, V) &=
\left((\emptyset,T)\rightarrow (S,\emptyset) \right) \in \cB_{\sigma,\tau}, \\
W&=\left((\emptyset,T)\rightarrow (S,\emptyset ) \right)_R \in
\mathcal{C}^{(\mu,\nu)}_{(\sigma,\tau)},
\end{split}
\end{equation*}
for some $\sigma,\tau \in\cP$. By Proposition \ref{tensor
decomposition}, there exist unique $\widetilde{U}\in \cB_{\sigma}$
and $\widetilde{V}\in\cB_{\tau}^\vee$ such that
$\widetilde{V}\otimes \widetilde{U}\equiv (U,V)$. By Proposition \ref{LR decomp}, we have
\begin{equation*}
W \longleftrightarrow (W_1,W_2)\in {\bf LR}^{\mu}_{\sigma
\lambda}\times {\bf LR}^{\nu}_{\tau\lambda }
\end{equation*}
for some $\lambda\in\cP$. By Proposition \ref{skew LR}, there exists
unique $X\in \cB_{\mu/\lambda}$ and $Y\in \cB_{\nu/\lambda}$ such
that
\begin{equation*}
\begin{split}
&\jmath(X)=\widetilde{U}, \ \ \ \jmath(X)_R=W_1, \\
&{\jmath}(Y)^\vee=\widetilde{V}, \ \ \ {\jmath}(Y)_R=W_2.
\end{split}
\end{equation*}
Now, we define
\begin{equation}
\psi_{\mu,\nu}(S\otimes T) =Y^\vee \otimes  X \in
\cB_{\nu/\lambda}^\vee\otimes \cB_{\mu/\lambda}.
\end{equation}
By construction, $\psi_{\mu,\nu}$ is bijective and commutes with
$\widetilde{x}_i$ for $x=e,f$ and $i\geq 1$. Hence we have the
following.

\begin{prop}\label{isomorphism-1+infty}
For $\mu,\nu\in\cP$, the map
\begin{equation*}
\psi_{\mu,\nu} : \cB_\mu\otimes \cB_\nu^\vee\longrightarrow
\bigsqcup_{\lambda\subset \mu,\nu}\cB_{\nu/\lambda}^\vee\otimes
\cB_{\mu/\lambda}
\end{equation*}
is an isomorphism of $\gl_{>0}$-crystals.
\end{prop}

\begin{ex}\label{insertion ex-3}{\rm Let
$S$ and
$T$ be the tableaux in Example \ref{insertion ex-1}.
Let
\begin{equation*}
X=
\begin{array}{ccc}
\bullet & \bullet   & 1  \\
\bullet & 1
\end{array}, \ \ \
Y=
\begin{array}{ccc}
\bullet & \bullet & 1   \\
\bullet & 2 \\
4
\end{array}.
\end{equation*}
Following the above notations, we have
\begin{equation*}
\begin{split}
& H_{(2,1)}\ast X=
\begin{array}{ccc}
{\bf 1} & {\bf 1}   & 1  \\
{\bf 2} & 1
\end{array} \ \stackrel{\text switching}{\leftrightsquigarrow} \
\begin{array}{ccc}
1 & 1   & {\bf 1}  \\
{\bf 1} & {\bf 2}
\end{array}= j(X)\ast j(X)_R= \widetilde{U}\ast W_1, \\
& H_{(2,1)}\ast Y =
\begin{array}{ccc}
{\bf 1} & {\bf 1} & 1   \\
{\bf 2} & 2 \\
4
\end{array} \ \stackrel{\text switching}{\leftrightsquigarrow} \
\begin{array}{ccc}
1 & 2 & {\bf 1}   \\
4 & {\bf 2} \\
{\bf 1}
\end{array}= {j}(Y)\ast  {j}(Y)_R=\widetilde{V}\ast
W_2,
\end{split}
\end{equation*}
where $\widetilde{U}$,
$\widetilde{V}$, $W_i$ ($i=1,2$) are as
in Examples \ref{insertion ex-1} and \ref{insertion ex-2}. Hence,
\begin{equation*}
\begin{split}
\psi_{\mu,\nu}(S\otimes T)&= Y^\vee \otimes X \\
&= \left(
\begin{array}{ccc}
\bullet & \bullet & 1   \\
\bullet & 2 \\
4
\end{array}
\right)^\vee \otimes \begin{array}{ccc}
\bullet & \bullet   & 1  \\
\bullet & 1
\end{array} \\
&=\begin{array}{ccc}
    &  &  4^\vee  \\
   & 2^\vee & \!\!\! \bullet  \\
 1^\vee   & \!\!\!\bullet  & \!\!\!\bullet
\end{array}
\otimes
\begin{array}{ccc}
\bullet & \bullet   & 1  \\
\bullet & 1
\end{array}.\\
\end{split}
\end{equation*}
}
\end{ex}

For a skew Young diagram and $\lambda/\mu$ and $k\geq 1$, we define
\begin{equation}\label{shifting}
\kappa_k :  \cB_{\lambda/\mu} \longrightarrow
\cB_{(\lambda+(1^k))/(\mu+(1^k))}
\end{equation}
by $\kappa_k(S)=S'$ with
\begin{equation*}
S'(i,j)=
\begin{cases}
S(i,j), & \text{if $i> k$}, \\
S(i,j-1), & \text{if $i\leq k$}.
\end{cases}
\end{equation*}

\begin{ex}{\rm
\begin{equation*}
\begin{split}
\kappa_1\left(
\begin{array}{ccc}
\bullet & \bullet & 1   \\
\bullet & 2 \\
1
\end{array}\right)
\ \ = \ \
\begin{array}{cccc}
\bullet & \bullet & \bullet & 1   \\
\bullet & 2 \\
1
\end{array} \ \ \ \ \
\kappa_2\left(
\begin{array}{ccc}
\bullet & \bullet & 1   \\
\bullet & 2 \\
1
\end{array}\right)
\ \ = \ \
\begin{array}{cccc}
\bullet & \bullet & \bullet & 1   \\
\bullet & \bullet & 2 \\
1
\end{array}
\end{split}
\end{equation*}
}
\end{ex}

For $k\geq 1$ and $\lambda\in\cP$, we put
\begin{equation*}
\begin{split}
\omega_{k}&=\epsilon_{1}+\cdots+\epsilon_{k}, \\
\omega_{\lambda}&=\lambda_1\epsilon_{1}+\lambda_{2}\epsilon_{2}+\cdots.
\end{split}
\end{equation*}
Now, we have the following combinatorial interpretation of the
embedding (\ref{embedding}) in terms of {\it sliding} skew tableaux
horizontally. It will play a crucial role in proving our main
theorem.
\begin{prop}\label{isomorphism-2+infty} For $\mu,\nu\in \cP$ and $k\geq 1$,
we have the following commutative diagram of $\gl_{>0}$-crystal
morphisms
$$
\begin{CD}
 \cB_{\mu}\otimes \cB_{\nu}^\vee @>\iota_{\omega_\mu,\omega_\nu}^{\omega_k} >> \cB_{\mu+(1^k)}\otimes
 \cB_{\nu+(1^k)}^\vee \\
 @V\psi_{\mu,\nu}VV @VV\psi_{\mu+(1^k),\nu+(1^k)}V \\
 \bigsqcup_{\lambda}\cB_{\nu/\lambda}^\vee\otimes
\cB_{\mu/\lambda} @>\kappa_k^\vee\otimes \kappa_k >>
 \bigsqcup_{\eta}\cB_{(\nu+(1^k))/\eta}^\vee\otimes
\cB_{(\mu+(1^k))/\eta}
\end{CD}
$$
where $\iota_{\omega_\mu,\omega_\nu}^{\omega_k}$ is the canonical
embedding in $($\ref{embedding}$)$ and
$\kappa^\vee_k=\vee\circ\kappa_k\circ\vee$.
\end{prop}
\pf Let $S\otimes T\in \cB_\mu\otimes \cB_\nu^\vee$ be given. We
keep the previous notations. Note that
\begin{equation*}
S\otimes u_{\omega_k}=S\otimes H_{(1^k)}\equiv S\{k\}:=\left(1\cdots
k \rightarrow\ S\right) \in \cB_{\mu+(1^k)},
\end{equation*}
\begin{equation*}
u_{-\omega_k}\otimes T=H_{(1^k)}^\vee\otimes T \equiv
T\{k\}:=\left(1\cdots k \rightarrow\ T^\vee\right)^\vee  \in
\cB^\vee_{\nu+(1^k)}.
\end{equation*}
Hence by (\ref{embedding}) we have
$\iota_{\omega_\mu,\omega_\nu}^{\omega_k}(S\otimes T)=S\{k\}\otimes
T\{k\}$. Since $S\{k\}\otimes T\{k\}\equiv S\otimes T$, we have
\begin{equation*}
(U\{k\},V\{k\}):=\left((\emptyset,T\{k\})\rightarrow
(S\{k\},\emptyset )\right) \equiv \left((\emptyset,T)\rightarrow
(S,\emptyset )\right)=(U, V),
\end{equation*}
which implies that $(U\{k\},V\{k\})=(U,V)$ by \cite[Lemma 5.1]{K09}.
Put
\begin{equation*}
W\{k\}=\left((\emptyset,T\{k\})\rightarrow (S\{k\},\emptyset
)\right)_R,
\end{equation*}
and suppose that
\begin{equation*}
W\{k\} \longleftrightarrow (W_1\{k\}, W_2\{k\})\in {\bf
LR}^{\mu+(1^k)}_{\sigma \lambda+(1^k)}\times {\bf
LR}^{\nu+(1^k)}_{\tau\lambda+(1^k) }.
\end{equation*}

Since $W$ is invariant under $\te_i$ and $\tf_i$ ($i\geq 1$), we may
assume that $(U,V)=(H^{>0}_\sigma, E^{>0}_\tau(n))$ for a
sufficiently large $n>k$ (see (\ref{Extremal weight vector})). As a
$\gl_{[s]}$-crystal element, $(U,V)$ is a highest weight element,
and $\sigma_n^{p}(U,V)=H^{>0}_{\zeta}$, where $p\geq \tau_1$ and
$\zeta=\sigma+(p-\tau_{n-i+1})_{i\geq 1}$ (see \cite[Section
4.1]{K09} for the definition of the map $\sigma_n$). This also
implies that $S=H^{>0}_{\mu}$. By \cite[Lemma 7.6]{Str}, $W\{k\}$ is
obtained from
\begin{equation}
\sigma_n^{-p}\Bigl[\bigl(\sigma_n^{p}(\emptyset,T\{k\})\rightarrow
(S\{k\},\emptyset)\bigr)_R\Bigr]
\end{equation}
by $180^\circ$-rotation and ignoring $\vee$'s in the entries. Since
$S\{k\}=H^{>0}_{\mu+(1^k)}$, we have
$\bigl(\sigma_n^{p}(\emptyset,T\{k\})\rightarrow
(S\{k\},\emptyset)\bigr)_R=\sigma_n^{p}(\emptyset,T\{k\})$. Now, it
is straightforward to check that
\begin{equation*}
W\{k\}=
\begin{array}{c}
  1 \\
  \vdots \\
  k
\end{array}
\ast \kappa_k(W).
\end{equation*}
This implies that
\begin{equation*}
\begin{split}
W_1\{k\} &= W_1 \ast \Sigma_k,\\
W_2\{k\} &= W_2 \ast \Sigma_k',
\end{split}
\end{equation*}
where $\Sigma_k$ and $\Sigma'_k$ are vertical strips of shape
$(\mu+(1^k))/\mu$ and $(\nu+(1^k))/\nu$ filled with $1,\ldots,k$,
respectively. Now, we have
\begin{equation*}
\begin{split}
\widetilde{U}\ast W_1\{k\}
=\widetilde{U}\ast W_1 \ast \Sigma_k & \rightsquigarrow
H_{\lambda}\ast X \ast \Sigma_k  \ \ \ \ \ \ \ \   \text{(switching $\widetilde{U}$ and $W_1$)}  \\
& \rightsquigarrow H_{\lambda+(1^k)}\ast \kappa_k(X) \ \ \ \text{(switching $X$ and $\Sigma_k$)},\\
\widetilde{V}\ast W_2\{k\}
=\widetilde{V}\ast W_2 \ast \Sigma'_k & \rightsquigarrow
H_{\lambda}\ast Y \ast \Sigma'_k  \ \ \ \ \ \  \ \ \text{(switching $\widetilde{V}$ and $W_2$)}  \\
& \rightsquigarrow H_{\lambda+(1^k)}\ast \kappa_k(Y) \ \ \ \text{(switching $Y$ and $\Sigma'_k$)}.
\end{split}
\end{equation*}
Therefore, it follows that
\begin{equation*}
\begin{split}
\psi_{\mu+(1^k),\nu+(1^k)}(\iota_{\omega_\mu,\omega_\nu}^{\omega_k}(S\otimes
T))
&=\psi_{\mu+(1^k),\nu+(1^k)}(S\{k\}\otimes T\{k\}))\\
&=\kappa_k(Y)^\vee\otimes \kappa_k(X) \\
&=\kappa_k^\vee\otimes \kappa_k\left( \psi_{\mu,\nu}(S\otimes
T)\right).
\end{split}
\end{equation*}
\qed

\subsection{}\label{modified crystal} Let $\cM$ be the set of
$\mathbb{N}\times \mathbb{N}$ matrices $A=(a_{ij})$
such that $a_{ij}\in\Z_{\geq 0}$ and $\sum_{i,j\geq 1}a_{i
j}<\infty$. Let $A=(a_{ij})\in \cM$ be given. For $i\geq 1$, the
$i$-th row $A_i=(a_{ij})_{j\geq 1}$ is naturally identified with a
unique semistandard tableau in $\cB_{(m_i)}$, where $m_i=\sum_{j\geq
1}a_{ij}$  and ${\rm wt}\left(A_i\right)=\sum_{j\geq
1}a_{ij}\epsilon_j$. Hence $A$ can be viewed as an element in
$\cB_{(m_1)}\otimes\ldots \otimes \cB_{(m_r)}$ for some $r\geq 0$.
This defines a $\gl_{>0}$-crystal structure on $\cM$. Now, we put
\begin{equation}
\widetilde{\cM}=\cM^\vee\times \cM,
\end{equation}
which can be viewed as a tensor product of $\gl_{>0}$-crystals.
Let $\mathcal{P}=\bigoplus_{i\geq 1}\Z\epsilon_i$ be the integral weight
lattice for $\gl_{>0}$. For $\omega\in \mathcal{P}$, let
\begin{equation*}
\widetilde{\cM}_\omega=\{\,(M^\vee, N) \in \widetilde{\cM}
\,|\,{\rm wt}(N^t) - {\rm wt}(M^t) =\omega \,\}.
\end{equation*}
Here $A^t$ denotes the transpose of $A\in\cM$. Then
$\widetilde{\cM}_\omega$ is a subcrystal of $\widetilde{\cM}$. Now,
we can state the main result in this section.

\begin{thm}\label{Momega} For $\omega\in \mathcal{P}$, we have
$$\widetilde{\cM}_\omega \simeq \B(\infty)\otimes T_{\omega}\otimes
\B(-\infty).$$
\end{thm}
\pf Let  $\mu,\nu\in\cP$ be such that
$\omega=\omega_\mu-\omega_\nu$. Suppose that
$\psi_{\mu,\nu}(S\otimes T)=Y^\vee \otimes X$ for $S\otimes T\in
\cB_{\mu}\otimes \cB_{\nu}^\vee$, where $\psi_{\mu,\nu}$ is the
isomorphism in Proposition \ref{isomorphism-1+infty}. Let
$M=(m_{ij})$ (resp. $N=(n_{ij})$) be the unique matrix in $\cM$ such
that the $i$-th row of $M$ (resp. $N$) is $\gl_{>0}$-equivalent to
the $i$-th row of $Y$ (resp. $X$). Since $\sum_{j\geq 1}m_{ij}$
(resp. $\sum_{j\geq 1}n_{ij}$) is equal to $x_i$ (resp. $y_i$) the
number of dots or boxes in the $i$-th row of $X$ (resp. $Y$) for $i\geq 1$
and $\omega=\sum_{i\geq 1}(x_i-y_i)\epsilon_i$  by Proposition
\ref{isomorphism-1+infty}, we have ${\rm wt}(N^t) - {\rm wt}(M^t)
=\omega$. Then we define
\begin{equation*}
\iota'_{\mu,\nu} : \cB_{\mu}\otimes \cB_{\nu}^\vee \longrightarrow
\widetilde{\cM}_\omega
\end{equation*}
by $\iota'_{\mu,\nu}(S\otimes T)=(M^\vee, N)$. By Proposition
\ref{isomorphism-1+infty}, it is easy to see that $\iota'_{\mu,\nu}$ is an
embedding and
\begin{equation*}\label{Mchar-1+infty}
\widetilde{\cM}_\omega=\bigcup_{\substack{\mu,\nu\in \cP \\
\omega_\mu-\omega_\nu=\omega}}{\rm Im}\iota'_{\mu,\nu}.
\end{equation*}
For $k\geq 1$, we have
$\iota'_{\mu,\nu}=\iota'_{\mu+(1^k),\nu+(1^k)}\circ
\iota^{\omega_k}_{\omega_\mu,\omega_\nu}$ by Proposition
\ref{isomorphism-2+infty}. Using induction, we have
\begin{equation*}\label{Mchar-2}
\iota'_{\mu,\nu}=\iota'_{\mu+\xi,\nu+\xi}\circ
\iota^{\omega_\xi}_{\omega_\mu,\omega_\nu} \ \ \ (\xi\in \cP).
\end{equation*}
Therefore, by (\ref{characterization of Btilde}), it follows that
$\widetilde{\cM}_\omega\simeq \B(\infty)\otimes T_{\omega}\otimes
\B(-\infty)$. \qed\vskip 5mm

\begin{cor}\label{MBiso} As a $\gl_{>0}$-crystal, we have
$$\B(\widetilde{U}_q(\gl_{>0})) \simeq \widetilde{\cM}.$$
\end{cor}
\pf It follows from $\widetilde{\cM}=\bigsqcup_{\omega\in
\mathcal{P}}\widetilde{\cM}_\omega$. \qed

For $A\in\cM$ and $i\geq 1$, we also define
\begin{equation}
\begin{split}
\te_i^t A= \left(\te_i A^t\right)^t, \ \ \tf_i^t A = \left(\tf_i
A^t\right)^t.
\end{split}
\end{equation}
Then $\cM$ has another $\gl_{>0}$-crystal structure with respect to
$\te_i^t$, $\tf_i^t$ and ${\rm wt}^t$, where ${\rm wt}^t(A)={\rm
wt}(A^t)$. By \cite{DK}, $\cM$ is a $(\gl_{>0},\gl_{>0})$-bicrystal,
that is, $\te_i,\tf_i$ on $\cM\cup\{{\bf 0}\}$ commute with
$\te_j^t, \tf_j^t$ for $i,j\geq 1$, and so is the tensor product
$\widetilde{\cM}=\cM^\vee\times \cM$. Now we have the following
Peter-Weyl type decomposition.
\begin{cor}\label{PeterWeyl+infty} As a
$(\gl_{>0},\gl_{>0})$-bicrystal, we have
\begin{equation*}
\B(\widetilde{U}_q(\gl_{>0}))\simeq\bigsqcup_{\mu,\nu\in\cP}\cB_{\mu,\nu}
\times \cB_{\mu,\nu}.
\end{equation*}
\end{cor}
\pf Note that the usual RSK correspondence gives an isomorphism of
$(\gl_{>0},\gl_{>0})$-bicrystals $\cM \simeq
\bigsqcup_{\lambda\in\cP}\cB_\lambda\times \cB_\lambda$. We assume
that $\te_i,\tf_i$ act on the first component, and $\te_j^t,
\tf_j^t$ act on the second component. The decomposition of
$\B(\widetilde{U}_q(\gl_{>0}))$ follows from Proposition \ref{tensor
decomposition}.\qed

\section{Extremal weight crystals of type $A_{\infty}$}
In this section, we describe the tensor product of
$\gl_\infty$-crystals $\B(\Lambda)\otimes \B(-\Lambda')$ for
$\Lambda, \Lambda'\in P^+$ in terms of extremal weight crystals.
\subsection{}
For a skew Young diagram $\lambda/\mu$, we put
\begin{equation}\label{Bmu}
\B_{\lambda/\mu}=SST_{\mathbb{Z}}(\lambda/\mu),
\end{equation}
and we identify $\B_{\lambda/\mu}^\vee$ with
$SST_{\mathbb{Z}^\vee}(\left(\lambda/\mu\right)^\vee)$. Note that
for $\mu\in\cP$, $\B_{\mu}$ has neither highest nor lowest weight
vector. It is shown in \cite{K09-2} that for $\mu,\nu\in\cP$,
$\B_{\mu}\otimes\B_{\nu}^\vee$ is connected,
$\B_{\mu}\otimes\B_{\nu}^\vee\simeq \B_{\nu}^\vee\otimes \B_{\mu}$,
and $\B_{\mu}\otimes\B_{\nu}^\vee \simeq
\B_{\sigma}\otimes\B_{\tau}^\vee$ if and only if
$(\mu,\nu)=(\sigma,\tau)$. Put
\begin{equation}
\B_{\mu,\nu}=\B_\mu\otimes \B_\nu^\vee.
\end{equation}
Note that $\B_{\mu,\nu}$ can be viewed as a limit of
$\B^{>r}_{\mu,\nu}$ ($r\rightarrow -\infty$) since
$\B^{>r}_{\mu,\nu}\simeq \left(\B_{\nu}^{>r}\right)^\vee\otimes
\B_{\mu}^{>r}$.

Let
$\Z_+^n=\{\,\lambda=(\lambda_1,\ldots,\lambda_n)\in\Z^n\,|\,\lambda_1\geq\cdots\geq\lambda_n\,\}$
be the set of generalized partitions of length $n$. For
$\lambda\in\Z_+^n$, we put
\begin{equation*}
\begin{split}
\Lambda_{\lambda}&
=\Lambda_{\lambda_1}+\cdots+\Lambda_{\lambda_n}\in P^+_n.
\end{split}
\end{equation*}

\begin{thm}[Theorem 4.6 in \cite{K09-2}]\label{extremal}
For $\Lambda\in P_n$ $(n\geq 0)$, there exist unique $\lambda\in
\Z_n^+$ and $\mu,\nu\in\cP$ such that
$$\B(\Lambda)\simeq \B_{\mu,\nu}\otimes\B(\Lambda_\lambda).$$
Here we assume that $\Lambda_\lambda=0$ when $n=0$.
\end{thm}

Note that $\{\B_{\mu,\nu}\otimes\B(\Lambda)\,|\,\Lambda\in P^+, \
\mu,\nu\in\cP\,\}$  forms a complete list of extremal weight
crystals of non-negative level up to isomorphism \cite[Proposition
3.12]{K09-2}.

\subsection{}\label{matrices}
For intervals $I,J$ in $\Z$, let $M_{I,J}$ be the set of $I\times J$
matrices $A=(a_{ij})$ with $a_{ij}\in\{\,0,1\,\}$. We denote by
$A_i$ the $i$-th row of $A$ for $i\in I$.

Suppose that $A\in M_{I,J}$ is given. For $k\in J^\circ$ and $i\in I$, we define
\begin{equation}
\tf_k A_i=
\begin{cases}
A_i-E_{ik}+E_{i k+1}, & \text{if $(a_{ik},a_{i\, k+1})=(1,0)$},\\
{\bf 0}, & \text{otherwise}.
\end{cases}
\end{equation}
Then we can regard $A_i$ as an element of a regular
$\gl_{\{k,k+1\}}$-crystal with highest weight
$\omega\in\{\,0,\epsilon_k,\epsilon_k+\epsilon_{k+1}\,\}$. Consider
the sequence $(\varepsilon_k(A_i))_{i\in I}$. We say that {\it $A$
is $k$-admissible} if there exist $L,L'\in I$ such that (1)
$\varepsilon_k(A_i)\neq 1$ for all $i<L$, and (2)
$\varepsilon_k(A_i)\neq -1$ for all $i>L'$. Note that if $I$ is
finite, then $A$ is $k$-admissible for all $k\in J$. Suppose that
$A$ is $k$-admissible. Then we can define $\tf_kA$ by regarding $A$
as $\overrightarrow{\bigotimes}_{i\in I}A_i $ and applying tensor
product rule of crystal, where the index $i$ in the tensor product
is increasing from left to right.

Let $\rho : M_{I,J} \longrightarrow M_{-J,I}$ be a bijection given
by $\rho(A)=(a_{-j,i})\in M_{-J,I}$, where $-J=\{\,-j\,|\,j\in
J\,\}$. For $l\in I^\circ$, we say that $A$ is {\it $l$-admissible}
if $\rho(A)$ is $l$-admissible in the above sense. If $A$ is
$l$-admissible, then we define
\begin{equation}
\begin{split}
\tE_l(A)=\rho^{-1}\left(\te_{l}\,\rho(A) \right), \ \ \tF_l(A)=\rho^{-1}\left(\tf_{l}\, \rho(A) \right)
\end{split}.
\end{equation}
If $A$ is both $k$-admissible and
$l$-admissible for some $l\in I^\circ$ and $k\in J^\circ$, then
\begin{equation}\label{commmuting}
\widetilde{x}_k\widetilde{X}_lA=\widetilde{X}_l\widetilde{x}_kA,
\end{equation}
where $x=e,f$ and $X=E,F$ \cite[Lemma 3.1]{K09-2}.

For convenience, let us say that {\it $A$ is $J$-admissible} (resp.
{\it $I$-admissible}) if $A$ is $k$-admissible (resp.
$l$-admissible) for all $k\in J^\circ$ (resp. $l\in I^\circ$).
Suppose that $A$ is $J$-admissible and $l$-admissible for some $l\in
I^\circ$. Then both $A$ and $\widetilde{X}_lA$ generate the same
$J$-colored  oriented graphs with respect to $\te_k$ and $\tf_k$ for
$k\in J^\circ$  whenever $\widetilde{X}_lA\neq {\bf 0}$ $(X=E,F)$
\cite[Lemma 3.2]{K09-2}. A similar fact holds when $A$ is
$I$-admissible and $k$-admissible for some $k\in J^\circ$.

If $I$ and $J$ are finite, then $M_{I,J}$ is a
$(\gl_I,\gl_J)$-bicrystal,  where the $\gl_I$-weight (resp.
$\gl_J$-weight) of $A=(a_{ij})\in M_{I,J}$ is given by $\sum_{i\in
I}\sum_{j\in J}a_{ij}\epsilon_i$ (resp. $\sum_{j\in J}\sum_{i\in
I}a_{ij}\epsilon_j$).

\subsection{}\label{Fock space}
For $n\geq 1$, let $\E^n$ be the subset of $M_{[n],\Z}$ consisting
of matrices $A=(a_{ij})$ such that $\sum_{i,j}a_{ij}<\infty$. It is
clear that $A$ is $\Z$-admissible for $A\in \E^n$. If we define
${\rm wt}(A)=\sum_{j\in \Z} \left(\sum_{i\in[n]} a_{ij}
\right)\epsilon_j$, then $\E^n$ is a regular $\gl_\infty$-crystal
with respect to $\te_k, \tf_k$ ($k\in\Z$) and ${\rm wt}$. For
$r\in\Z$ and $\lambda\in \cP$ with $\lambda_1\leq n$, let
$A^\circ_{\lambda}(r)=(a^\circ_{ij})\in\E^n$ and
$A^\diamond_{\lambda}(r)=(a^\diamond_{ij})\in\E^n$ be such that for
$i\in [n]$ and $j\in\Z$
\begin{equation}\label{Alambda positive level}
\begin{split}
&a^\circ_{ij}=1 \ \ \Longleftrightarrow  \ \ 1+r \leq j\leq
\lambda'_{n-i+1}+r,\\
&a^\diamond_{ij}=1 \ \ \Longleftrightarrow  \ \ r- \lambda'_{i}+1 \leq j\leq
r.
\end{split}
\end{equation}
Then $C(A^\ast_{\lambda}(r))\simeq \B_\lambda$ ($\ast = \circ, \diamond$).

For $n\geq 1$, let $\F^n$ be the set of matrices $A=(a_{ij})$ in
$M_{[n],\Z}$ such that for each $i\in [n]$, $a_{ij}=1$ if $j\ll 0$
and $a_{ij}=0$ if $j\gg 0$. Note that $A$ is $\Z$-admissible for
$A\in \F^n$. If we define ${\rm
wt}(A)=n\Lambda_0+\sum_{j>0}\sum_{i\in
[n]}a_{ij}\epsilon_j+\sum_{j\leq 0}\sum_{i\in
[n]}(a_{ij}-1)\epsilon_j$, then $\F^n$ is a regular
$\gl_\infty$-crystal with respect to $\te_k, \tf_k$ ($k\in\Z$) and
${\rm wt}$. For $\lambda\in\Z_+^n$, let
$A^\circ_{\lambda}=(a^\circ_{ij})\in \F^n$ and
$A^\diamond_\lambda=(a^\diamond_{ij})\in \F^n$ be such that for
$i\in [n]$ and $j\in\Z$
\begin{equation}\label{Alambda}
\begin{split}
&a^\circ_{ij}=1 \ \ \Longleftrightarrow \ \ j\leq \lambda_{n-i+1}, \\
&a^\diamond_{ij}=1 \ \ \Longleftrightarrow \ \ j\leq \lambda_{i}.
\end{split}
\end{equation}
Then  $C(A_\lambda^\ast)\simeq \B(\Lambda_\lambda)$
($\ast=\circ,\diamond)$.

On the other hand, for $A=(a_{ij})\in\E^n$ or $\F^n$, $A$ is
$[n]$-admissible.  Hence, $\widetilde{E}_l$ and $\widetilde{F}_l$
($l=1,\ldots, n-1$) commute with $\te_k$ and $\tf_k$ ($k\in \Z$).

For $A=(a_{ij})\in\E^n$ or $\F^n$, we will identify its dual crystal
element $A^\vee$   with the matrix $(1-a_{ij})$ since $A^\vee$ and
$(1-a_{ij})$ generate the same $\Z$-colored graph.\vskip 2mm

\subsection{}
Let $m,n$ be non-negative integers with $m\geq n$. In the rest of
this section, we fix $\mu\in\Z_+^m$ and $\nu\in \Z_+^n$. We assume
that $\B(\Lambda_\mu)=C(A_\mu^\circ)\subset \F^m$,
$\B(-\Lambda_\nu)=C((A_\nu^\diamond)^\vee)\subset
\left(\F^n\right)^\vee$ and hence
$$\B(\Lambda_\mu)\otimes \B(-\Lambda_\nu)\subset \F^m\otimes
\left(\F^n\right)^\vee\subset M_{[m+n],\Z}.$$
By
\cite[Proposition 4.5]{K09-2}, $\F^m\otimes \left(\F^n\right)^\vee$
is a disjoint union of extremal weight $\gl_\infty$-crystals of
level $m-n$,  and hence so is $\B(\Lambda_\mu)\otimes
\B(-\Lambda_\nu)$.

For $r\in\Z$, we define $\B^{> r}(\mu,\nu)$ to be the set of
$A=(a_{ij})\in  M_{[m+n],\Z}$ such that
\begin{equation*}
a_{ij}=
\begin{cases}
1, & \text{for $i\in [n]$ and $j\leq r$}, \\
0,  & \text{for $i\in m+[n]$ and $j\leq r$}.
\end{cases}
\end{equation*}
We have
\begin{equation*}
\begin{split}
&\B^{> r+1}(\mu,\nu)\subset \B^{> r}(\mu,\nu),\\
&\B(\Lambda_\mu)\otimes \B(-\Lambda_\nu)=\bigcup_{r\in\Z}\B^{> r}(\mu,\nu).
\end{split}
\end{equation*}

Let $r$ be such that $r<\min\{\mu_m,\nu_n\}$ so that
$\mu-(r^m)=(\mu_i-r)_{1\leq i\leq m}$ and
$\nu-(r^n)=(\nu_i-r)_{1\leq i\leq m}$ are partitions. Then $\B^{>
r}(\mu,\nu)\neq \emptyset$ and as a $\gl_{>r}$-crystal,
\begin{equation}\label{tensor product >r}
\B^{> r}(\mu,\nu)\simeq \B_{(\mu-(r^m))'}^{>r}\otimes
\left(\B_{(\nu-(r^n))'}^{>r}\right)^\vee.
\end{equation}
\vskip 2mm

Now, let $A\in \B^{> r}(\mu,\nu)$ be given and
$C^{>r}(A)$ the connected component in $\B^{>r}(\mu,\nu)$
including $A$ as a $\gl_{>r}$-crystal. By Proposition \ref{LR rule
for >r}, we have
\begin{equation*}
C^{>r}(A)\simeq \B_{\sigma,\tau}^{>r}
\end{equation*}
for some $\sigma,\tau\in \cP$ with $\sigma_1\leq m$ and $\tau_1\leq
n$.  Let $C(A)$ be the connected component in $\B(\Lambda_\mu)\otimes \B(-\Lambda_\nu)$
including $A$ as a $\gl_{\infty}$-crystal. Then
\begin{equation*}
C(A)\simeq \B_{\zeta,\eta}\otimes\B(\Lambda_\xi)
\end{equation*}
for some $\zeta,\eta\in \cP$ and $\xi\in\Z_+^{m-n}$ by Theorem \ref{extremal}.

\begin{lem}\label{weight correspondence} Under the above hypothesis, we have
\begin{equation*}
\begin{split}
\zeta&=(\sigma'_{m-n+1},\ldots,\sigma'_m)', \\
\eta&= \tau, \\
\xi&= (\sigma'_1,\ldots,\sigma'_{m-n})+(r^{m-n}).
\end{split}
\end{equation*}
\end{lem}
\pf Let $A$ be as above. For intervals $I,J \subset\Z$,  let
$A_{I,J}$ denote the $I\times J$-submatrix of $A$. Choose $s\gg r$
so that
\begin{equation*}
a_{ij}=
\begin{cases}
0, & \text{if $i\in [m]$ and $j>s$}, \\
1, & \text{if $i\in m+[n]$ and $j>s$.}
\end{cases}
\end{equation*}
Note that $A$ is $[m+n]$-admissible. Considering $A_{[m+n],[r+1,s]}$ as an element
of a $(\gl_{[r+1,s]},\gl_{[m+n]})$-bicrystal, $A$ is connected to a unique matrix
$A'=(a'_{ij})$ satisfying
\begin{equation*}
\begin{cases}
 a'_{ij}=a_{ij}, & \text{for $i\in [m+n]$ and $j\not\in [r+1,s]$}, \\
  a'_{i-1\, j}=0, & \text{if $a'_{ij}=0$ for $i\neq 1$ and $j\in [r+1,s]$},\\
  a'_{i\, j+1}=0, & \text{if $a'_{ij}=0$ for $i\in [m+n]$ and $j+1\in [r+1,s]$}.
\end{cases}
\end{equation*}
Equivalently, $A'$ is a $\gl_{[r+1,s]}$-highest weight element and a $\gl_{[m+n]}$-lowest weight element.
Then we have
\begin{equation*}
\begin{split}
&C^{>r}(A'_{[m],\Z}) \simeq \B_\alpha^{>r}, \\
&C^{>r}(A'_{m+[n],\Z})\simeq \left(\B_\beta^{>r}\right)^\vee,
\end{split}
\end{equation*}
where $\alpha=(\alpha_k)_{k\geq 1}$ and $\beta=(\beta_k)_{k\geq
1}\in \cP$  are given by $\alpha_k=\sum_{i=1}^m a'_{i, r+k}$ for
$1\leq k\leq s-r$ and $\beta_k=\sum_{i=1}^n (1-a'_{m+i, s-k+1})$ for
$1\leq k\leq s-r$. Since $A'_{[m+n],[r+1,\infty)}$ is
$\gl_{>r}$-equivalent to $H^{>r}_\alpha\otimes E^{>r}_{\beta}(s)$
(see (\ref{Extremal weight vector})), we have $C^{>r}(A')\simeq
\B_{\alpha,\beta}^{>r}$ by Theorem \ref{extremal +infty} and hence
$(\alpha,\beta)=(\sigma,\tau)$ since $C^{>r}(A')\simeq
C^{>r}(A)\simeq \B_{\sigma,\tau}^{>r}$.

Let $A''=(a''_{ij})\in M_{[m+n],\Z}$ be such that
\begin{equation*}
\begin{split}
& A''_{[n],\Z}=A^\circ_{\zeta}(r)\in\E^n, \\
& A''_{n+[n],\Z}=\left(A^\diamond_{\eta}(s)\right)^\vee\in \left(\E^n\right)^\vee, \\
& A''_{2n+[m-n],\Z}=A^\circ_{\xi}\in\F^{m-n},
\end{split}
\end{equation*}
where $\zeta=(\sigma'_{m-n+1},\ldots,\sigma'_m)'$  and $\eta= \tau'$
and $\xi= (\sigma'_1,\ldots,\sigma'_{m-n})+(r^{m-n})$ (see (\ref{Alambda positive level}) and (\ref{Alambda})).
Then for $L\ll 0\ll L'$, we have
\begin{equation*}
S_{w'}S_{w}(A'_{[m+n],[L,L']})=A''_{[m+n],[L,L']},
\end{equation*}
where
\begin{equation*}
\begin{split}
w&=\left(r_{n-2}\cdots r_{1}\right)\cdots \left(r_{m+n-2}\cdots
r_{m-1}\right)\left(r_{m+n-1}\cdots r_m\right),\\
w'&=\left(r_{2n}\cdots r_{m+n}\right)\cdots \left(r_{n+2}\cdots
r_{m+2}\right)\left(r_{n+1}\cdots r_{m+1}\right),
\end{split}
\end{equation*}
and $S_w$, $S_{w'}$ are the corresponding operators on a regular
$\gl_{[m+n]}$-crystal $M_{[m+n],[L,L']}$  with respect to $\tE_i$, $\tF_i$'s. This
implies that $A'$ is $\gl_{[L,L']}$-equivalent to $A''$. Since $L$
and $L'$ are arbitrary, $A'$ is
$\gl_\infty$-equivalent to $A''$. Since
\begin{equation*}
\begin{split}
&C(A''_{[2n],\Z})\simeq \B_{\zeta,\eta}, \\
&C(A''_{2n+[m-n],\Z})\simeq \B(\Lambda_\xi),
\end{split}
\end{equation*}
we have
\begin{equation*}
C(A)\simeq C(A')\simeq C(A'')\simeq
\B_{\zeta,\eta}\otimes\B(\Lambda_\xi).
\end{equation*}
This completes the proof. \qed
\vskip 5mm

For $\zeta,\eta\in \cP$ and $\xi\in\Z_+^{m-n}$, let
$m^{(\mu,\nu)}_{(\zeta,\eta,\xi)}(r)$ be the number of connected
components $C$ in $\B(\Lambda_\mu)\otimes \B(-\Lambda_\nu)$ such that
\begin{itemize}
\item[(1)] $C\cap \B^{>r}(\mu,\nu)\neq \emptyset$,

\item[(2)] $C\simeq \B_{\zeta,\eta}\otimes\B(\Lambda_\xi)$.
\end{itemize}

\begin{cor}\label{LR identity} Under the above hypothesis,
\begin{itemize}
\item[(1)] if $\xi_{m-n}< r$, then $m^{(\mu,\nu)}_{(\zeta,\eta,\xi)}(r)=0$,

\item[(2)] if $\xi_{m-n}\geq r$, then
$$m^{(\mu,\nu)}_{(\zeta,\eta,\xi)}(r)=c^{((\mu-(r^m))',(\nu-(r^n))')}_{(\sigma,\eta)},$$ where
$\sigma=\left[\left(\xi-(r^{m-n})\right)\cup \zeta'\right]'$.
\end{itemize}
\end{cor}
\pf It follows from (\ref{tensor product >r}) and Lemma \ref{weight
correspondence}. \qed\vskip 2mm

The following lemma shows that $m^{(\mu,\nu)}_{(\zeta,\eta,\xi)}(r)$ stabilizes as $r\rightarrow -\infty$.
\begin{lem}\label{stability}
For $\zeta,\eta\in \cP$ and $\xi\in\Z_+^{m-n}$, there exists
$r_0\in\Z$ such that
\begin{equation*}
m^{(\mu,\nu)}_{(\zeta,\eta,\xi)}(r)=m^{(\mu,\nu)}_{(\zeta,\eta,\xi)}(r_0),
\end{equation*}
for $r\leq r_0$.
\end{lem}
\pf For $r\in\Z$ with $r\leq \min\{\mu_m,\nu_n\}$, put
\begin{equation*}
\mathcal{C}^{(\mu,\nu)}_{(\zeta,\eta,\xi)}(r)
=\bigsqcup_{\lambda\in\cP}{\rm\bf LR}^{(\mu-(r^m))'}_{\sigma
\lambda}\times {\rm\bf LR}^{(\nu-(r^n))'}_{\eta\, \lambda},
\end{equation*}
where $\sigma=\left[\left(\xi-(r^{m-n})\right)\cup \zeta'\right]'$.
Then
\begin{equation*}
\mathcal{C}^{(\mu,\nu)}_{(\zeta,\eta,\xi)}(r-1)
=\bigsqcup_{\lambda\in\cP}{\rm\bf LR}^{(\mu-(r^m))'\cup\{m\}}_{\overline{\sigma}\,\,
\lambda\cup\{n\}}\times {\rm\bf LR}^{(\nu-(r^n))'\cup\{n\}}_{\eta\,\, \lambda\cup\{n\}},
\end{equation*}
where $\overline{\sigma}=\left[\left(\xi-(r^{m-n})+(1^{m-n})\right)\cup
\zeta'\right]'$.

By Proposition \ref{LR rule for >r} and Corollary \ref{LR identity}, we have
\begin{equation*}
\left|\mathcal{C}^{(\mu,\nu)}_{(\zeta,\eta,\xi)}(r)\right|=
c^{((\mu-(r^m))',(\nu-(r^n))')}_{(\sigma,\eta)}(r)=m^{(\mu,\nu)}_{(\zeta,\eta,\xi)}(r).
\end{equation*}
For a sufficiently small $r$, we define a map
\begin{equation*}
\theta_r : \mathcal{C}^{(\mu,\nu)}_{(\zeta,\eta,\xi)}(r)
\longrightarrow \mathcal{C}^{(\mu,\nu)}_{(\zeta,\eta,\xi)}(r-1).
\end{equation*}
as follows;

\noindent \textsc{Step 1.} Suppose that $S_1\in {\rm \bf
LR}^{(\mu-(r^m))'}_{\sigma \lambda}$ is given. Put $\ell=\xi_{m-n}-r$.

Define $T_1$ to be the tableau in ${\rm\bf
LR}^{(\mu-(r^m))'\cup\{m\}}_{\overline{\sigma}\ \lambda\cup\{n\}}$, which
is obtained from $S_1$ as follows;
\begin{itemize}
\item[(1)] The entries of $T_1$ in the $i$-th row ($1\leq i\leq \ell$)
 is equal to those in $S_1$.
\item[(2)] The entries of $T_1$ in the $(\ell+1)$-st row is
given by
$$a_1+1\leq a_2+1\leq \cdots\leq a_n+1,$$
where $a_1\leq a_2\leq\cdots\leq a_n$ are the entries in the
$\ell$th row in $S_1$.

\item[(3)] Let $S'_1$ (resp. $T'_1$) be the subtableau of $S_1$ (resp.
$T_1$) consisting of its $i$ th row for $\ell<i$ (resp. $\ell+1<i$). Then we define
\begin{equation*}
T'_1(p+1,q)=
\begin{cases}
S'_1(p,q), & \text{if $S'_1(p,q)\leq a_1$},\\
S'_1(p,q)+1, & \text{if $S'_1(p,q) > a_1$},\\
\end{cases}
\end{equation*}
for $(p,q)$ in the shape of $S_1$.
\end{itemize}

\noindent \textsc{Step 2.} Let $S_2\in {\rm\bf
LR}^{(\nu-(r^n))'}_{\eta \lambda}$ be given. Applying the same
argument as in Step 1 (when $m=n$). we obtain $T_2\in {\rm\bf
LR}^{(\nu-(r^n))'\cup\{n\}}_{\eta (\lambda'+(1^n))'}$. Now we define
\begin{equation*}
\theta_r(S_1,S_2)=(T_1,T_2)\in
\mathcal{C}^{(\mu,\nu)}_{(\zeta,\eta,\xi)}(r-1).
\end{equation*}

By definition of $\theta_r$, it is not difficult to see that
$\theta_r$ is one-to-one. Also, we observe that for $\lambda\in\cP$
\begin{equation*}
{\rm\bf LR}^{(\mu-(r^m))'}_{\sigma \lambda}\times {\rm\bf
LR}^{(\nu-(r^n))'}_{\eta \lambda} \neq \emptyset \ \
\Longleftrightarrow  \ \ {\rm\bf LR}^{(\mu-(r^m))'\cup\{m\}}_{\overline{\sigma}
\,\, \lambda\cup\{n\}}\times {\rm\bf LR}^{(\nu-(r^n))'\cup\{n\}}_{\eta
\,\, \lambda\cup\{n\}} \neq \emptyset.
\end{equation*}
If $r$ is sufficiently small, then we have $(n)\subset \lambda$  for
$\lambda\in\cP$ such that ${\rm\bf
LR}^{(\mu-(r^m))'\cup\{m\}}_{\overline{\sigma} \lambda}\times
{\rm\bf LR}^{(\nu-(r^n))'\cup\{n\}}_{\eta \lambda} \neq \emptyset$,
which implies that $\theta_r$ is onto.  Therefore, $\theta_r$ is a
bijection and $m^{(\mu,\nu)}_{(\zeta,\eta,\xi)}(r)$ stabilizes as
$r\rightarrow -\infty$. \qed \vskip 3mm

\begin{thm}\label{main theorem} Suppose that $m\geq n$.
For $\mu\in\Z_+^m$ and $\nu\in\Z_+^n$, we have
\begin{equation*}
\B(\Lambda_\mu)\otimes \B(-\Lambda_\nu) \simeq
\bigsqcup_{\substack{\zeta,\eta\in\cP \\ \zeta_1,\eta_1\leq
n}}\bigsqcup_{\xi\in\Z_+^{m-n}}
\B_{\zeta,\eta}\otimes\B(\Lambda_\xi)^{\oplus
m^{(\mu,\nu)}_{(\zeta,\eta,\xi)}}
\end{equation*}
with
\begin{equation*}
m^{(\mu,\nu)}_{(\zeta,\eta,\xi)}=
\sum_{\lambda\in\cP}c^{\mu+(k^m)}_{ \sigma\lambda
}c^{\nu+(k^n)}_{\eta\,\lambda},
\end{equation*}
where $k$ is a sufficiently large integer and
$\sigma=\left(\xi+(k^{m-n})\right)\cup \zeta'$.
\end{thm}
\pf For $\zeta,\eta\in \cP$ and $\xi\in\Z_+^{m-n}$, let
$m^{(\mu,\nu)}_{(\zeta,\eta,\xi)}$ be the number of connected
components in $\B(\Lambda_\mu)\otimes \B(-\Lambda_\nu)$ isomorphic to
$\B_{\zeta,\eta}\otimes\B(\Lambda_\xi)$. Then by
Lemma \ref{stability}, we have
\begin{equation*}
m^{(\mu,\nu)}_{(\zeta,\eta,\xi)}=m^{(\mu,\nu)}_{(\zeta,\eta,\xi)}(r)
\end{equation*}
for some $r\in\Z$. By Corollary \ref{LR rule for >r} and Corollary \ref{LR identity}, we have
\begin{equation*}
m^{(\mu,\nu)}_{(\zeta,\eta,\xi)}=
\sum_{\lambda\in\cP}c^{\mu+(k^m)}_{ \sigma\lambda
}c^{\nu+(k^n)}_{\eta\lambda},
\end{equation*}
where $k=-r$ and  $\sigma=\left(\xi+(k^{m-n})\right)\cup
\zeta$. \qed

The decomposition when $m\leq n$ can be obtained by taking the dual
crystal of the decomposition in Theorem \ref{main theorem}.

\section{Combinatorial description of the Level zero part of $\B(\widetilde{U}_q(\gl_\infty))$}

\subsection{}\label{n-n tensor bijection}
For $\mu,\nu\in\Z_+^n$ ($n\geq 1$), let us describe the
decomposition of $\B(\Lambda_\mu)\otimes \B(-\Lambda_\nu)$ in a
bijective way. We assume that $\B(\Lambda_\mu)=C(A_\mu^\circ)\subset
\F^m$, $\B(-\Lambda_\nu)=C((A_\nu^\diamond)^\vee)\subset
\left(\F^n\right)^\vee$

Suppose that $A\in \B(\Lambda_\mu)$ and $A'\in \B(-\Lambda_\nu)$ are
given. Choose $r\in\Z$ such that $A\otimes A'\in \B^{>r}(\mu,\nu)$.
Let $S^{>r}\otimes T^{>r}\in \B_{(\mu-(r^m))'}^{>r}\otimes
\left(\B_{(\nu-(r^n))'}^{>r}\right)^\vee$ correspond to $A\otimes
A'$ under (\ref{tensor product >r}). Note that the set of entries in
the $i$-th column of $S^{>r}$ (from the right) is
$\{\,j\,|\,a_{ij}=1,\ j>r\,\}$, and the set of entries in the $i$-th
column of $T^{>r}$ (from the right) is $\{\,j^\vee\,|\,a_{ij}=0,\
j>r\,\}$. Now we define
\begin{equation}
\begin{split}
&\psi^\infty_{\mu,\nu}(A\otimes A')
=\psi^{>r}_{(\mu-(r^m))',(\nu-(r^n))'}(S^{>r}\otimes T^{>r}),
\end{split}
\end{equation}
where $\psi^{>r}_{(\mu-(r^m))',(\nu-(r^n))'}$ denotes the map in
Proposition \ref{isomorphism-1+infty} corresponding to
$\gl_{>r}$-crystals.

\begin{prop}\label{isomorphism-1}
For $\mu,\nu\in\Z_+^n$, the map
\begin{equation*}
\psi^\infty_{\mu,\nu} : \B(\Lambda_\mu)\otimes
\B(-\Lambda_\nu)\longrightarrow \bigsqcup_{\alpha,\beta}
\B_{\alpha}^\vee\otimes\B_{\beta}
\end{equation*}
is an isomorphism of $\gl_\infty$-crystals, where the union is over
all skew Young diagrams $\alpha$ and $\beta$ such that
$\alpha=\left(\nu-(r^n)\right)'/\lambda$ and
$\beta=\left(\mu-(r^n)\right)'/\lambda$ for some $r\leq
\min\{\mu_n,\nu_n\}$ and $\lambda\in\cP$.
\end{prop}
\pf  It suffices to show that $\psi^\infty_{\mu,\nu}(A\otimes A')$ does
not depend on the choice of $r$.  Keeping the above notations, we have
\begin{equation*}
\begin{split}
(U^{>r}, V^{>r}) &=
\left((\emptyset,T^{>r})\rightarrow (S^{>r},\emptyset ) \right) \in \B^{>r}_{\sigma,\tau}, \\
W^{>r}&=\left((\emptyset,T^{>r})\rightarrow (S^{>r},\emptyset )
\right)_R \in
\mathcal{C}^{(\left(\mu-(r^n)\right)',\left(\nu-(r^n)\right)')}_{(\sigma,\tau)},
\end{split}
\end{equation*}
for some $\sigma,\tau \in\cP$.
By Proposition \ref{tensor
decomposition}, there exist unique $\widetilde{U}^{>r}\in
\B^{>r}_{\sigma}$ and
$\widetilde{V}^{>r}\in\left(\B^{>r}_{\tau}\right)^\vee$ such that
$\widetilde{V}^{>r}\otimes \widetilde{U}^{>r}\equiv (U^{>r},
V^{>r})$, and by Proposition \ref{LR decomp}
\begin{equation*}
W^{>r} {\longleftrightarrow} (W^{>r}_1,W^{>r}_2)\in {\bf
LR}^{\left(\mu-(r^n)\right)'}_{\sigma \lambda}\times {\bf
LR}^{\left(\nu-(r^n)\right)'}_{\tau\lambda }
\end{equation*}
for some $\lambda\in\cP$. Then by definition of
$\psi^{>r}_{(\mu-(r^n))',(\nu-(r^n))'}$, we have
\begin{equation*}
\psi^\infty_{\mu,\nu}(A\otimes A') =Y^\vee \otimes  X \in
\B_{\left(\nu-(r^n)\right)'/\lambda}^\vee\otimes \B_{\left(\mu-(r^n)\right)'/\lambda},
\end{equation*}
where
\begin{equation*}
\begin{split}
&\jmath(X)=\widetilde{U}^{>r}, \ \ \ \jmath(X)_R=W_1^{>r}, \\
&{\jmath}(Y)^\vee=\widetilde{V}^{>r}, \ \ \ {\jmath}(Y)_R=W_2^{>r}.
\end{split}
\end{equation*}

Now, suppose that $$S^{>r-1}\otimes T^{>r-1}\in
\B_{(\mu-(r^n))'\cup\{n\}}^{>r}\otimes
\left(\B_{(\nu-(r^n))'\cup\{n\}}^{>r}\right)^\vee$$ is
$\gl_{>r-1}$-equivalent to $A\otimes A'$. Then
$$S^{>r-1}=(\underbrace{r-1\cdots r-1}_{n})\ast S^{>r},\ \
T^{>r-1}=T^{>r}\ast (\underbrace{(r-1)^\vee\cdots
(r-1)^\vee}_{n}),$$ and
$$\left((\emptyset,T^{>r-1})\rightarrow (S^{>r-1},\emptyset ) \right)
=\left((\emptyset,T^{>r})\rightarrow (S^{>r},\emptyset )
\right)=(U^{>r}, V^{>r}).$$  This implies that $(U^{>r-1},
V^{>r-1})=(U^{>r}, V^{>r})$.

Suppose that $W^{>r}=W^{>r}_+\ast W^{>r}_-$, where $W^{>r}_+$ (resp.
$W^{>r}_-$) is the subtableau of $W^{>r}$ consisting of positive
(resp. negative) entries. By definition of the insertion, it is
straightforward to check that
\begin{itemize}
\item[(1)] $W^{>r-1}_-=W^{>r}_-$,

\item[(2)] $W^{>r-1}_+=(\underbrace{\sigma'_n+1 \cdots \sigma'_1+1}_{n}) \ast W^{>r}_+[1]$,
\end{itemize}
where $W^{>r}_+[1]$ is the tableau obtained from $W^{>r}_+$ by
increasing each entry by $1$. Since
$\imath(W^{>r-1}_+)=W^{>r-1}_{1}$, we have
$$W^{>r-1}_1=\Sigma_n\ast W_1^{>r}[1],$$
where $\Sigma_n$ is the horizontal strip of shape
$\sigma\cup\{n\}/\sigma$ filled with $1$, and $W_1^{>r}[1]$ is the
tableau obtained from $W_1^{>r}$ by increasing each entry by $1$.
Here, we assume that the shape of $W_1^{>r}$ is
$(\mu-(r^n))'\cup\{n\}/\sigma\cup\{n\}$. Now, we have
\begin{equation*}
\begin{split}
\widetilde{U}^{>r}\ast W^{>r-1}_1 &=\widetilde{U}^{>r}\ast \Sigma_n \ast W^{>r}_{1}[1] \\
& \rightsquigarrow(\underbrace{1 \cdots 1}_{n}) \ast \widetilde{U}^{>r}\ast W^{>r}_{1}[1]
\ \ \ \ \ \text{(switching $\widetilde{U}^{>r}$ and $\Sigma_n$)}\\
& \rightsquigarrow (\underbrace{1 \cdots 1}_{n}) \ast H_\lambda[1] \ast X
\ \ \ \ \ \ \ \  \ \ \ \text{(switching $\widetilde{U}^{>r}$ and $W^{>r}_{1}[1]$)}\\
&= H_{\lambda\cup\{n\}}\ast X.
\end{split}
\end{equation*}
This implies that $X$ does not depend on $r$. Similarly, we have
$$W^{>r-1}_2=\Sigma'_n\ast W^{>r}_2[1],$$
where $\Sigma'_n$ is the horizontal strip of shape
$\tau\cup\{n\}/\tau$ filled with $1$, and
\begin{equation*}
\begin{split}
\widetilde{V}^{>r}\ast W^{>r-1}_2 &=\widetilde{V}^{>r}\ast \Sigma'_n \ast W^{>r}_{2}[1] \\
& \rightsquigarrow(\underbrace{1 \cdots 1}_{n}) \ast \widetilde{V}^{>r}\ast W^{>r}_{2}[1]
 \ \ \ \ \ \text{(switching $\widetilde{V}^{>r}$ and $\Sigma'_n$)} \\
& \rightsquigarrow (\underbrace{1 \cdots 1}_{n}) \ast H_\lambda[1] \ast Y
\ \ \ \ \ \ \ \  \ \ \ \text{(switching $\widetilde{V}^{>r}$ and $W^{>r}_{2}[1]$)} \\
&=H_{\lambda\cup\{n\}}\ast Y.
\end{split}
\end{equation*}
This also implies that $Y$ does not depend on $r$. Therefore,
$\psi^\infty_{\mu,\nu}$ is well-defined.

Since $\psi^\infty_{\mu,\nu}$ is one-to-one and commutes with
$\te_i$ and $\tf_i$ ($i\in \Z$) by construction, it is an
isomorphism of $\gl_\infty$-crystals. \qed

\begin{ex}{\rm
Let $\mu=(2,2,1)$ and $\nu=(3,2,1)$. Consider
\begin{equation*}
\begin{split}
A&=
\begin{array}{cccccccccccc}
  & \!\!-3 & \!\!-2 & \!\!-1 & \ 0 & \ 1 & \ 2 & \ 3 & \ 4 & \ 5 & \   &   \\
  \cdots & \bullet & \bullet & \bullet & \bullet & \cdot & \bullet & \cdot & \cdot & \cdot &  \cdots & 1\\
  \cdots & \bullet & \bullet & \bullet & \bullet & \bullet & \cdot & \bullet & \cdot & \cdot &  \cdots & 2 \\
  \cdots & \bullet & \bullet & \bullet & \bullet & \bullet & \bullet & \cdot & \cdot & \cdot &  \cdots & 3
\end{array}\ \in \B(\Lambda_\mu)\subset \F^3,\\
A'&=
\begin{array}{cccccccccccc}
  & \!\!-3 & \!\!-2 & \!\!-1 & \ 0 & \ 1 & \ 2 & \ 3 & \ 4 & \ 5 & \   &   \\
  \cdots & \cdot & \cdot & \cdot & \cdot & \cdot & \cdot & \bullet & \cdot & \bullet &  \cdots & 1\\
  \cdots & \cdot & \cdot & \cdot & \cdot & \bullet & \cdot & \cdot & \bullet & \bullet &  \cdots & 2 \\
  \cdots & \cdot & \cdot & \cdot & \cdot & \bullet & \cdot & \bullet & \bullet & \bullet &  \cdots & 3
\end{array}\ \in \B(-\Lambda_\nu)\subset \left(\F^3\right)^\vee,
\end{split}
\end{equation*}
where $\bullet$ and $\cdot$ denote $1$ and $0$ in matrix,
respectively. Then $A\otimes A'\in \B^{>0}(\mu,\nu)$. Suppose that
as a $\gl_{>0}$-crystal element $A$ (resp. $A'$) is equivalent to
$S^{>0}$ and $T^{>0}$. Then $S^{>0}=S$ and $T^{>0}=T$, where $S$ and
$T$ are tableaux in Example \ref{insertion ex-1}. Hence, by Example
\ref{insertion ex-3} we have
\begin{equation*}
\begin{split}
\psi^\infty_{\mu,\nu}(A\otimes A')&=
\begin{array}{ccc}
    &  &  4^\vee  \\
   & 2^\vee & \!\!\! \bullet  \\
 1^\vee   & \!\!\!\bullet  & \!\!\!\bullet
\end{array}
\otimes
\begin{array}{ccc}
\bullet & \bullet   & 1  \\
\bullet & 1
\end{array}.\\
\end{split}
\end{equation*}
}
\end{ex}

\subsection{}
Let us give an explicit description of $\B(\infty)\otimes T_\Lambda
\otimes \B(-\infty)$ for $\Lambda\in P_0$. For this, we define an
analogue of (\ref{shifting}) for $\gl_\infty$-crystals. Suppose that
$\mu\in\Z_+^n$ is given. For $k\in\Z$, let $\mu\cup\{k\}$ be the
generalized partition in $\Z_+^{n+1}$ given by rearranging
$\mu_1,\ldots,\mu_n$ and $k$. For $r\leq \mu_n$, we assume that the
columns in $\left(\mu-(r^n)\right)'\in \cP$ are enumerated from the
left, and the row indices are enumerated by $r+1, r+2,\ldots$ from
the top. For a skew Young diagram
$\alpha=\left(\mu-(r^n)\right)'/\lambda$ and and $S\in \B_{\alpha}$,
we also denote by $S(i,j)$ the entry in $S$ located in the $i$-th
row and the $j$-th column. For $k\in\Z$, we define $\kappa_k :
SST_{\Z}(\alpha) \rightarrow SST_{\Z}(\kappa_k(\alpha))$, where
$\kappa_k(\alpha)=\left((\mu\cup\{k\})-(r^{n+1})\right)'/\,(\lambda+(1^{k-r}))$
and $\kappa_k(S)=S'$ is given by $S'(i,j)=S(i,j)$ if $i> k$, and
$S(i,j-1)$ if $i\leq k$. We put
$\kappa_k^\vee=\vee\circ\kappa_k\circ\vee$. If $k<r$, then we assume
that $\alpha=\left(\mu-(s^n)\right)'/\lambda+(n^{r-s})$ for $s\leq
k$.

By applying the argument in Proposition \ref{isomorphism-2+infty} to
Proposition \ref{isomorphism-1} with a little modification, we
obtain the following.

\begin{prop}\label{isomorphism-2} For $\mu,\nu\in \Z_+^n$ $(n\geq 1)$ and $k\in\Z$,
we have the following commutative diagram of $\gl_\infty$-crystal
morphisms.
$$
\begin{CD}
 \B(\Lambda_\mu)\otimes \B(-\Lambda_\nu) @>\iota^{\Lambda_k}_{\Lambda_\mu,\Lambda_\nu}>>
 \B(\Lambda_\mu+\Lambda_k)\otimes
 \B(-\Lambda_k-\Lambda_\nu) \\
 @V\psi^\infty_{\mu,\nu}VV @VV\psi^\infty_{\mu\cup\{k\},\nu\cup\{k\}}V \\
 \bigsqcup_{\alpha,\beta} \B_{\alpha}^\vee\otimes\B_{\beta} @>\kappa_k^\vee\otimes \kappa_k >>
 \bigsqcup_{\gamma,\delta} \B_{\gamma}^\vee\otimes\B_{\delta}
\end{CD}
$$
\end{prop}\vskip 3mm

Let $\M$ be the set of $\Z\times \Z$ matrices $A=(a_{ij})$ such that
$a_{ij}\in\Z_{\geq 0}$ and $\sum_{i,j\in\Z}a_{i j}<\infty$. Let
$A=(a_{ij})\in \M$ be given. As in Section \ref{modified crystal},
we have a $(\gl_\infty,\gl_\infty)$-bicrystal structure on $\M$ with
respect to $\te_i, \tf_i$ and $\te_j^t, \tf_j^t$ for $i,j\in \Z$.
Now, we put
\begin{equation}
\begin{split}
\widetilde{\M}&=\M^\vee\times \M,\\
\widetilde{\M}_\Lambda&=\{\,(M^\vee, N) \in \widetilde{\M}
\,|\,{\rm wt}(N^t) - {\rm wt}(M^t) =\Lambda \,\} \ \ \ \ (\Lambda\in
P_0).
\end{split}
\end{equation}
Note that $\widetilde{\M}$ can be viewed as a tensor product of
$(\gl_\infty,\gl_\infty)$-bicrystals and $\widetilde{\M}_\Lambda$ is
a subcrystal of $\widetilde{\M}$ with respect to $\te_i, \tf_i$. By
Proposition \ref{isomorphism-2}, we have the following combinatorial
realization, which is our second main result. The proof is almost
the same as in Theorem \ref{Momega}.

\begin{thm} For $\Lambda\in P_0$, we have
$$\widetilde{\M}_\Lambda \simeq \B(\infty)\otimes T_{\Lambda}\otimes \B(-\infty).$$
\end{thm}

Let $\B(\widetilde{U}_q(\gl_\infty))_0=\bigsqcup_{\Lambda\in
P_0}\B(\infty)\otimes T_{\Lambda}\otimes \B(-\infty)$ be the level
zero part of $\B(\widetilde{U}_q(\gl_\infty))$. Since
$\widetilde{\M}=\bigsqcup_{\Lambda\in P_0} \widetilde{\M}_\Lambda$
and ${\M}\simeq \bigsqcup_{\lambda\in\cP}\B_\lambda\times
\B_\lambda$ as a $(\gl_\infty,\gl_\infty)$-bicrystal, we obtain the
following immediately.
\begin{cor}\label{MBiso} As a $\gl_\infty$-crystal, we have
$$\B(\widetilde{U}_q(\gl_\infty))_0\simeq \widetilde{\M}.$$
\end{cor}

\begin{cor} As a
$(\gl_\infty,\gl_\infty)$-bicrystal, we have
\begin{equation*}
\B(\widetilde{U}_q(\gl_\infty))_0\simeq\bigsqcup_{\mu,\nu\in\cP}\B_{\mu,\nu}\times\B_{\mu,\nu}.
\end{equation*}
\end{cor}

In \cite{BN}, Beck and Nakajima proved a Peter-Weyl type
decomposition of the level zero part of
$\B(\widetilde{U}_q(\frak{g}))$ for a quantum affine algebra
$\frak{g}$ of finite rank, where the bicrystal structure is given by
star crystal structure, say $\te_i^\ast$ and $\tf_i^\ast$, induced
from the involution on $\widetilde{U}_q(\frak{g})$, usually denoted
by $\ast$ \cite{Kas94'}. Based on some computation, we give the
following conjecture.

\begin{conj}
The crystal structure on $\B(\widetilde{U}_q(\gl_{>0}))$ and
$\B(\widetilde{U}_q(\gl_\infty))_0$ of type $A_{+\infty}$ and
$A_\infty$ with respect to $\te_i^t$ and $\tf_i^t$ is compatible
with the dual of the $\ast$-crystal strucutre with respect to
$\te_i^\ast$ and $\tf_i^\ast$. That is, $\te_i^t=\tf_i^\ast$ and
$\tf_i^t=\te_i^\ast$ for all $i$.
\end{conj}

{\small

\end{document}